\documentclass[a4paper,times new roman, 12pt]{amsart}

\usepackage{multirow}
\usepackage{amsfonts}
\usepackage{graphicx}
\usepackage{amsmath}
\usepackage{amssymb}
\usepackage[matrix,arrow,curve]{xy}
\usepackage[colorlinks=true, allcolors=blue]{hyperref}
\usepackage[usenames]{color}
\usepackage[mathmode, boxsize=3.7em]{ytableau}
\usepackage[square,sort,comma,numbers]{natbib}
\usepackage[utf8]{inputenc}
\usepackage[T1]{fontenc}
\newtheorem{thm}{theorem}[section]
\newtheorem{theorem}[thm]{Theorem}

\newtheorem{lemma}[thm]{Lemma}
\newtheorem{corollary}[thm]{Corollary}

\begin{document}
	
\title[Polynomial Identities of Certain Metabelian Leibniz Algebras]
{Polynomial Identities and Codimensions of Certain Three-Dimensional Metabelian Non-Lie Leibniz Algebras}

\author[Fertunani]{Luis Fertunani}
\thanks{L. Fertunani was partially supported by a Ph.D. fellowship from CNPq, Brazil, and by FAEPEX under Grant No. 2478/25.}
\address{Department of Mathematics, Universidade Estadual de Campinas (UNICAMP), 13083-859 Campinas, SP,  Brazil}
\email{l245102@dac.unicamp.br}

\author[Fideles]{Claudemir Fideles}
\thanks{C. Fidelis was supported by FAPESP grant No.~2025/26906-3}
\address{Department of Mathematics, Universidade Estadual de Campinas (UNICAMP), 13083-859 Campinas, SP,  Brazil}
\email{fideles@unicamp.br}

\author[Muniz]{Airton Muniz}
\thanks{A. Muniz was partially supported by a Ph.D. fellowship from CNPq, Brazil, and by FAEPEX under Grant No. 2478/25.}
\address{Department of Mathematics, Universidade Estadual de Campinas (UNICAMP), 13083-859 Campinas, SP,  Brazil}
\email{a217012@dac.unicamp.br}

	\begin{abstract}

Over an arbitrary field, we conduct an extensive study of the polynomial identities and codimensions of three-dimensional metabelian non-Lie Leibniz algebras. We investigate five one-parameter families and seven isolated algebras, determining finite bases for their corresponding \(T\)-ideals and providing explicit bases for the associated relatively free algebras. Due to the unique structure of each algebra considered, their explicit $T$-ideals and bases were obtained through tailored computational techniques of independent interest. While these results cover a broad class of three-dimensional metabelian non-Lie Leibniz algebras, a complete treatment over an arbitrary field remains dependent on the still open classification of such algebras. As a consequence of our analysis, we show that polynomial identities do not determine the isomorphism class in this setting: there exist non-isomorphic three-dimensional metabelian non-Lie Leibniz algebras satisfying exactly the same polynomial identities.

		\medskip
		
		\medskip
		
	\end{abstract}
	
	\maketitle
	
\section{Introduction}

Leibniz algebras were first seen in the papers of A.M. Blokh \cite{Blo65} and they were called right $D$-algebras, highlighting their close relationship with derivations. The concept of Leibniz algebras was introduced by Jean-Louis Loday \cite{loday1993version} in 1993 as a natural non-antisymmetric generalization of Lie algebras. While Lie algebras require the bracket operation to be antisymmetric and satisfy the Jacobi identity, Leibniz algebras relax the antisymmetry condition but retain a form of the Leibniz rule, or derivation property, expressed by the Leibniz identity. This generalization was motivated by the search for algebraic structures that arise naturally in homological algebra and noncommutative geometry, expanding the landscape of algebraic systems beyond Lie algebras. Loday’s paper laid the foundations for the theory, establishing the basic definitions and initial properties of Leibniz algebras. Additionally, the identification of Leibniz algebra cohomology with an Ext or Tor functor related to the Leibniz enveloping algebra, for more details see \cite{loday1993universal,zhang2015cohomology}.

The classification of low-dimensional algebras has been extensively studied in the literature. For further background, we refer the interested reader to \cite{bekbaev2021complete,SBI,C,bai2003automorphisms,kaygorodov2025polynomial} and the references therein, although many other relevant sources are available. In addition, subalgebras, idempotents, ideals, and quasi-units of two-dimensional algebras are investigated in \cite{ahmed2020subalgebras}.

In this context, the study of low-dimensional algebras plays a fundamental role in the development of the theory. Finite dimensional Leibniz algebra not only highlights the novel phenomena that emerge from the lack of antisymmetry, but also provides a foundational framework for understanding higher-dimensional cases. For example, in \cite{doi:10.1142/S0219199717500043}, the authors demonstrate how such can be associated with combinatorial structures of dimension 2.

On another front of research, it is well known that many important algebras are
\emph{PI-algebras}, that is, algebras satisfying a nontrivial set of polynomial
identities. Several classes of nonassociative algebras---including Leibniz, Jordan,
and Novikov algebras---are fundamental examples of PI-algebras. For instance, an
algebra $A$ is called a \emph{Leibniz algebra} if it satisfies the identity
\begin{equation}\label{ideLeibniz}
    x(yz) = (xy)z - (xz)y, \qquad \forall x, y, z \in A.
\end{equation}

It is important to note that two conventions appear in the literature.
An algebra satisfying \eqref{ideLeibniz} is called a \emph{right Leibniz algebra},
whereas a \emph{left Leibniz algebra} satisfies
\[
    (xy)z = x(yz) - y(xz), \qquad \forall x, y, z \in A.
\]
Throughout this paper, we adopt the first convention. Therefore, all algebras considered here are right Leibniz algebras in the standard terminology.

In recent years, many authors have established results relating algebraic structures and their polynomial identities. One natural question, called isomorphism problem, is: Are algebras that satisfy the same identities isomorphic?

Of course, any two isomorphic algebras satisfy the same identities. In general, however, the converse does not hold. For example, the algebras $M_{2}(\mathbb{R})$ (the algebra of $2 \times 2$ real matrices) and $\mathbb{H}_{\mathbb{R}}$ (the real quaternion algebra) satisfy the same identities but are not isomorphic. On the other hand, if the ground field $K$ is algebraically closed, it is well known that two finite-dimensional simple algebras satisfy the same identities if and only if they are isomorphic (see \cite{shestakov2011polynomial,bahturin2019distinguishing}). In dimension two, this equivalence holds in certain situations, particularly for Lie,  Leibniz and Jordan algebras. For a detailed discussion of this phenomenon in the context of Jordan algebras, we refer the reader to \cite{diniz2023two}.

Despite the extensive research in this area, very little is known about the explicit form of the identities satisfied by a given algebra. The monographs \cite{bahturin2021identical, giambruno2005polynomial}, along with the references therein, provide a thorough overview of the results obtained thus far. However, it should be emphasized that the list of algebras for which the polynomial identities are explicitly known remains very short and easy to reproduce. In the case of Lie algebras, known results include the identities of $UT_n(K)$ (as a Lie algebra), and also of $sl_2(K)$ over infinite fields \cite{razmyslov1994identities,vasilovskii1989basis}. Over a field of characteristic diﬀerent from two identities for $gl_2(K)$ (second order matrix algebra) coincide with identites for $sl_2(K)$, and a basis for this latter algebra over a finite field can be found in \cite{semenov1992basis}. These appear to be essentially the only well-understood cases. In the same context, Koshlukov \cite{koshlukov1997weak} determined a basis of identities for the natural representation $sl_2(K)\to M_2(K)$. Here we recall that weak identities can be defined in a similar way for other classes of algebras that are not ``very far'' from associative ones, and they play an important role in studying polynomial identities for Lie, alternative and Jordan algebras. In \cite{sheina1978metabelian},
Sheina determined the identities in certain finite metabelian Lie algebras. Finally, within the scope of studies on Leibniz algebras, we highlight the papers \cite{de2023polynomial,luis025polynomial}, which provides a treatment of the polynomial identities satisfied by their null-filiform version.

Given the difficulty of obtaining explicit descriptions of polynomial identities, it is natural to consider related invariants that may provide further information about them. Among these, codimensions play a particularly important role, as they provide a quantitative measure of polynomial identities and offer valuable insight into the structure of the corresponding algebras.

The notion of codimension for an algebra $A$ was introduced by Regev and leads to the study of the asymptotic behavior of its codimension sequence $(c_n(A))_{n \geq 1}$. For an associative PI-algebra, this sequence is exponentially bounded, and the limit
$$
\lim_{n \to \infty} \sqrt[n]{c_n(A)}
$$
exists and is a non-negative integer called the exponent of $A$. If $A$ is a finite-dimensional Lie algebra, then the codimension satisfies $c_n(A) \leq (\dim A)^{n+1}$. The codimensions of two-dimensional nonassociative algebras over a field of characteristic zero were studied in \cite{Giambruno20073405}. The authors proved that the codimension sequence of a two-dimensional nonassociative algebra $A$ is either bounded by $n + 1$ or grows exponentially as $2^n$. The restriction to fields of characteristic zero is essential, as the arguments rely on the representation theory of the symmetric group.

We now return to the setting of polynomial identities that will be the focus of this paper. Recall that an algebra satisfying the identity
$$
(x_{1}x_{2})(x_{3}x_{4})=0
$$
is called \emph{metabelian}; accordingly, throughout the paper we refer to this identity as the \emph{metabelian identity}.

Our main goal is to explicitly determine the polynomial identities of a broad class of three-dimensional metabelian non-Lie Leibniz algebras over an arbitrary field. More precisely, we consider five one-parameter families and seven isolated algebras, for which we determine finite bases of polynomial identities and explicit bases of the corresponding relatively free algebras. We also study their codimension sequences, proving that they are linearly bounded and that the corresponding PI-exponents are equal to~\(1\).

It is worth emphasizing that the determination of these finite bases and relatively free algebras is far from uniform. Due to the non-antisymmetric nature of Leibniz multiplication and the structural diversity across the parameterized families, a single standard algorithm or generic reduction does not apply. Instead, each algebra and family required a dedicated case-by-case analysis involving bespoke computational strategies and fine-grained reduction techniques to handle the specific interactions within their $T$-ideals. These technical methods hold independent interest for the computational study of nonassociative PI-algebras.

A complete classification, up to isomorphism, of three-dimensional non-Lie Leibniz algebras over arbitrary fields is beyond the scope of the present paper and, to the best of our knowledge, is not yet available in full generality. Some progress in this direction will be presented in the forthcoming work \cite{cordeiro2026classification}. Accordingly, our approach is guided primarily by polynomial identities rather than by isomorphism classes. Whenever necessary, distinct isomorphism classes can be identified by means of elementary invariants, such as the right annihilator and the dimensions of the terms of the lower central series.

%Beyond the explicit determination of polynomial identities, relatively free algebras, and codimension sequences for the classes considered here, our results reveal an interesting phenomenon: polynomial identities do not, in general, determine the isomorphism class of a three-dimensional metabelian non-Lie Leibniz algebra. Indeed, among the algebras studied in this paper, we exhibit non-isomorphic algebras satisfying exactly the same polynomial identities. Thus, even in dimension three, classification up to isomorphism and classification up to PI-equivalence lead to genuinely different problems.

Beyond providing explicit bases and invariants, our findings establish a remarkable phenomenon: \emph{polynomial identities do not determine the isomorphism class of three-dimensional metabelian non-Lie Leibniz algebras}. We exhibit explicit pairs of non-isomorphic algebras that share the exact same $T$-ideal. Thus, even in dimension three, classification up to isomorphism and classification up to PI-equivalence diverge significantly.

%From a broader perspective, this phenomenon places the present results within the general program of understanding varieties of algebras through their polynomial identities and numerical invariants. In PI-theory, codimension growth and the associated PI-exponent provide a bridge between explicit identities and the asymptotic structure of the variety generated by an algebra. In the Leibniz setting, where comparatively few explicit bases of identities are known, concrete descriptions of \(T\)-ideals and relatively free algebras therefore provide natural test cases for investigating how algebraic structure is reflected in codimension growth and PI-equivalence. In particular, the linear bounds obtained here identify a family of Leibniz varieties of polynomial growth and open the way to studying whether growth conditions can be used to organize broader classes of finite-dimensional Leibniz algebras.

%We expect that this point of view may contribute to a more systematic PI-theory for Leibniz algebras, in which explicit identities, relatively free objects, codimension growth, and PI-equivalence complement the usual classification by isomorphism. More generally, such an approach may be useful in comparing identity-theoretic invariants across different classes of nonassociative algebras and in understanding which structural properties are, or are not, detected by polynomial identities.

From a broader perspective, these results contribute to the overarching program of classifying varieties of nonassociative algebras via identity-theoretic and numerical invariants. In Leibniz theory, where explicit $T$-ideals are rare, our concrete descriptions establish a benchmark class of varieties with polynomial codimension growth, paving the way for future systematic investigations into how asymptotic invariants reflect core algebraic properties.

\section{Definitions and preliminary results} \label{Preliminaries}

In this paper, $K$ denotes an arbitrary field and any necessary restrictions will be stated explicitly. In particular, if $K$ is finite, we denote its cardinality by $|K|$. All algebras considered henceforth are taken over $K$. In this work, we are especially interested in Leibniz algebras, that is, algebras satisfying the Leibniz identity given in \eqref{ideLeibniz}.

Let $X = \{x_1, x_2, x_3, \dotsc\}$ be a countable infinite set of variables. We denote by $\mathcal{L}\langle X\rangle$ the free Leibniz algebra freely generated by $X$. Moreover, the elements of $\mathcal{L}\langle X\rangle$ will be called polynomials. Actually, by Leibniz identity \eqref{ideLeibniz}, any polynomial in $\mathcal{L}\langle X\rangle$ can be expressed as a
linear combination of left-normed monomials, i.e., monomials of the form
$(\ldots ((x_{i_1}x_{i_2})x_{i_3}) \ldots x_{i_{i-1}})$.
For this reason, we will omit the brackets in left-normed monomials. For example, the above monomial will be denoted simply by $x_{i_1}\cdots x_{i_m}$.

Given $f = f(x_1,\dots,x_n) \in \mathcal{L}\langle X\rangle$ and a Leibniz algebra $L$, an evaluation of $f$ in $L$ is a map $\varphi : \{x_1,\dots,x_n\} \to L$. We say that $f$ is a polynomial identity of $L$ if $f(a_1,\dots,a_n)=0$ for all $a_1,\dots,a_n \in L$, or equivalently, if for every evaluation $\varphi$ we have
\(
f(\varphi(x_1),\dots,\varphi(x_n)) = 0.
\)
The set of all polynomial identities of $L$, denoted by $T(L)$, is a $T$-ideal, that is, an ideal invariant under all endomorphisms of $\mathcal{L}\langle X\rangle$. The $T$-ideal generated by a subset $S \subset \mathcal{L}\langle X\rangle$ is the intersection of all $T$-ideals containing $S$, and it is usually denoted by $\langle S\rangle_T$. In this case, we say that $S$ is a basis of a $T$-ideal $I$ if $I = \langle S\rangle_T$. Moreover, $f$ is said to be a consequence of $S$ if $f$ belongs to the $T$-ideal generated by~$S$.

It is a classical result in PI-theory that, over arbitrary fields, every $T$-ideal is generated by its regular polynomials, i.e., polynomials in which all variables appear in every monomial, though not necessarily with the same degree. Moreover, the choice of generators depends on the ground field. For instance, in characteristic $0$, every $T$-ideal is generated by its multilinear polynomials. Over an infinite field of positive characteristic, one must instead consider multihomogeneous polynomials rather than multilinear ones.

Assume that the ground field $K$ is finite, with $|K|=q$. It is well known that the polynomials
\[
y_1y_2-y_2y_1 \quad \text{and} \quad y^q-y
\]
form a basis for the polynomial identities of $K$ viewed as an algebra. Moreover, the set of all monomials
\[
y_{i_1}^{\alpha_1}y_{i_2}^{\alpha_2}\cdots y_{i_t}^{\alpha_t},
\]
where $i_1<\cdots<i_t$ and $0\leq \alpha_j<q$ for all $j=1,\ldots,t$, forms a basis of the relatively free algebra $K[Y]/T(K)$ as a vector space over $K$. Here $K[Y]$ denotes the free associative algebra generated by the set $Y$.

The polynomial $f(x_{1},\dots, x_{n})\in \mathcal{L}\langle X\rangle$ is multilinear if each indeterminate $x_{1},\dots, x_{n}$ appears in every monomial of $f(x_{1},\dots, x_{n})$ exactly once. We denote by $P_n$ the subspace of $ \mathcal{L}\langle X\rangle$ of the multilinear polynomials in $n$ indeterminates. The $n$-th codimension of a Leibniz algebra is  
$$
c_n(L)=\dim \frac{P_n}{P_n\cap T(L)}.
$$
If the sequence $(c_n(L))_{n\geq 1}$ is exponentially bounded then one can consider the bounded sequence $\sqrt[n]{c_n(L)}$ with lower and upper limits
	$$
	\begin{array}{ccc} 
	\underline{\mbox{exp}}(L)= \liminf_{n \to \infty} \sqrt[n]{c_n(L)} & \mbox{and} & \overline{\mbox{exp}}(L)= \limsup_{n \to \infty} \sqrt[n]{c_n(L)}\\ 
	\end{array}
	$$
	called the lower and upper PI-exponents of $L$, respectively. If $\underline{\mbox{exp}}(L)=\overline{\mbox{exp}}(L)$ then the exponent of $L$ is $\mbox{exp}(L)= \lim_{n \to \infty} \sqrt[n]{c_n(L)}$.

\section{Identities, codimension and Relatively Free Algebras}

The classification of three-dimensional Leibniz algebras up to isomorphism already exhibits considerable structural diversity, even in such a low dimension. This diversity becomes particularly apparent from the viewpoint of polynomial identities, since different isomorphism classes may require substantially different techniques for the determination of their identities.

Over the complex field, three-dimensional Leibniz algebras were classified up to isomorphism by Ayupov and Omirov in \cite{ayupov19993}. This classification was recently revisited by Kaygorodov and Lopatin in \cite{kaygorodov2025polynomial}, who refined the original result and provided corrections to some of the cases therein. Further results and different approaches to the classification and structure of low-dimensional Leibniz algebras can be found in the monograph \cite{SBI} and in the papers \cite{rikhsiboev2012classification, kurdachenko2023structure, rakhimov2018algorithm}, among others.

A complete classification of three-dimensional Leibniz algebras over arbitrary fields falls outside the scope of the present paper. Such a classification has recently been obtained and will appear in the forthcoming work \cite{cordeiro2026classification}. The algebras investigated below are drawn from this classification and constitute the classes relevant to our study of polynomial identities.

For convenience, we now present the three-dimensional non-Lie Leibniz algebras that will be considered throughout this paper.
\[
\begin{array}{rllll}
L_1: & e_1 e_3 = -2 e_1, & e_2 e_2 = e_1, & e_3 e_2 = e_2, & e_2 e_3 = - e_2, \qquad (\operatorname{char}K\neq 2); \\
L_2: & e_1 e_3 = \alpha e_1, & e_3 e_2 = e_2, & e_2 e_3 = - e_2; \\
L_3: & e_3 e_3 = e_1, & e_3 e_2 = e_2, & e_2 e_3 = - e_2; \\
L_5: & e_2 e_2 = e_1, & e_3 e_3 = \alpha e_1; \\
L_6: & e_2 e_2 = e_1, & e_3 e_3 = \alpha e_1, & e_2 e_3 = e_1; \\
L_{7}: & e_2 e_3 = e_1; & \\
L_9: & e_1 e_3 = e_2, & e_2 e_3 = \alpha e_1; \\
L_{10}: & e_1 e_3 = e_2, & e_2 e_3 = \alpha e_1 + e_2; \\
L_{11}: & e_1 e_3 = e_1, & e_2 e_3 = e_2; \\
L_{12}: & e_3 e_3 = e_1, & e_1 e_3 = e_2; \\
L_{13}: & e_3 e_3 = e_1, & e_1 e_3 = e_1 + e_2; \\
N_1: & e_2 e_2 = e_1,  & e_3 e_3 = \alpha e_1, & e_3 e_2 = e_2, & e_2 e_3 = e_2, \qquad (\operatorname{char}K=2);\\
 N_3: &  e_3 e_3 = e_1;  & \\
 N_4: & e_1 e_3 = e_1,  & e_3 e_3 =  e_1. 
\end{array}
\]
Here $\{e_1,e_2,e_3\}$ denotes a basis of the algebra and $\alpha \in K$.

Firstly, if $K$ is a field of characteristic $2$, the algebra $L_1$ is replaced by $N_1$. Moreover, over quadratically closed fields, the algebras $L_5(\alpha)$ and $L_9(\alpha)$ are represented by $L_5(1)$ and $L_9(1)$, respectively. In the classification adopted in this paper, we write $L_4=L_5(1)$ and $L_8=L_9(1)$.

Note that the parameter in $L_2$ must be nonzero, otherwise, the algebra is Lie. Furthermore, for the algebras $L_5(\alpha)$, $L_6(\alpha)$,  $L_9(\alpha)$, and $L_{10}(\alpha)$, the case $\alpha=0$ yields algebras that are either isomorphic to previously listed algebras or trivial extensions. More precisely,
\[
L_6(0)\cong L_{7}, \qquad
L_9(0)\cong L_{7},
\]
while $L_5(0)$ is a trivial extension of a two-dimensional algebra studied in \cite{FERREIRADOSSANTOS20251} (denoted by $\textbf{T}_2$) and $L_{10}(0) \cong N_4$.
In this paper, we focus only on metabelian algebras. Accordingly, algebras of type $L_1$ and $N_1$ are excluded from our analysis.
Therefore, for the algebras that we will study throughout the paper, we can assume that $\alpha\neq 0$.

Notably, the Leibniz algebras $L_5$, $L_6$, $L_7$ and $L_{12}$ can be handled directly. These algebras are nilpotent: $L_5$, $L_6$ and $L_{7}$ have nilpotency index three, while $L_{12}$ is null-filiform. It is also important to note that $L_{5}$ is abelian. To summarize, we have the following result.

\begin{theorem}
Over any field, the identities of the algebras $L_6$ and $L_{7}$ coincide and are generated by the polynomial $x_1 x_2 x_3$, while 
\[T(L_{5}) = \langle x_1 x_2 - x_2 x_1, x_1 x_2 x_3 \rangle_T.\]
\end{theorem}

It is worth noting that $L_6$ and $L_{7}$ are not isomorphic, since their right annihilators have different dimensions. Thus, the previous theorem may be viewed as a counterexample to the isomorphism problem, that is, PI-equivalence does not imply isomorphism for Leibniz algebras of dimension greater than two.

Throughout the paper, the expression $x_1^{r_1} \cdots \widehat{x_j^{r_j}} \cdots x_m^{r_m}$ will denote the monomial obtained by omitting the factor $x_j^{r_j}$.

To simplify notation, we set
\[
ab^{(n)} := a\underbrace{bb\cdots b}_{n\ \text{times}},
\]
for all $a,b$ in a Leibniz algebra. This notation should not be confused with $ab^n$, which we interpret as
\[
ab^{n} := a(\underbrace{bb\cdots b}_{n\ \text{times}}).
\]
Thus, in $ab^{(n)}$ the product is left-normed, whereas in $ab^n$ the bracketing is right-normed due to the parentheses. This distinction is essential in nonassociative contexts.

  \begin{theorem}
Let $K$ be a field. Then the $T$-ideal of polynomial identities of the algebra $L_{12}$ is generated by:
\begin{enumerate}
\item[(i)] the polynomial $x_1 (x_2 x_3)$ and $x_1x_2x_3 - x_2x_1x_3$, if $|K| \geq 3$; and,
\item[(ii)] if $|K| = 2$, the polynomial $x_1 (x_2 x_3)$ and $x_1x_2x_3 - x_2x_1x_3$, together with the identity $x_1^{(2)} x_2 - x_1 x_2^{(2)}$.
\end{enumerate}
\end{theorem}

Finally, a direct computation shows that \(N_3\) and \(N_4\) are PI-equivalent to \(\mathbf{T}_2\) and \(\mathcal{L}_4\), respectively, where
\begin{equation}\label{eq:L4}
\mathcal{L}_4:\qquad e_1e_2=e_2e_2=e_1.
\end{equation}
Moreover, \(\mathcal{L}_4\), and hence \(N_4\), satisfies exactly the same polynomial identities as \(L_{11}\). Since the arguments are essentially the same, we postpone the proof to the corresponding discussion of \(L_{11}\).

We also record the polynomial identities of \(\mathbf{T}_2\). By arguments analogous to those used in \cite[Proposition~3.1]{FERREIRADOSSANTOS20251}, one obtains that
$$
[x_1,x_2]=x_1x_2-x_2x_1
\qquad\text{and}\qquad
(x_1x_2)x_3
$$
form a minimal basis for the \(T\)-ideal of polynomial identities of \(\mathbf{T}_2\) over an arbitrary field \(K\).

\medskip

\subsection{\texorpdfstring{The algebra $L_2$}{The algebra L2}}\

In the algebra \(L_2\), the multiplication on the basis \(\{e_1,e_2,e_3\}\) is given by  
\[
    e_1 e_3 = \alpha e_1, \qquad  e_3 e_2 = e_2, \qquad  e_2 e_3 = - e_2.
\]

Under this assumption, it is straightforward to verify that the polynomial \((x_1 x_2)(x_3 x_4)\) and the standard polynomial
\begin{equation}\label{s_2}
    s_3(x_1,x_2,x_3)
    = x_1x_2x_3 - x_1x_3x_2 - x_2x_1x_3 + x_2x_3x_1 + x_3x_1x_2 - x_3x_2x_1
\end{equation}
are identities for \(L_2\).

Now, if $|K| = q < \infty$, we add the following identities to the list:
\(
x^{q+1} - x^2,
\)
and the identities
\begin{align*}
&x_1 x_2^{(2)} x_3
- x_3 x_2^{(2)} x_1 x_3^{(q-1)}
+ x_3^{(q)} x_2^{(2)} x_1
- x_1 x_3^{(q)} x_2^{(2)}, \\[6pt]
&x_1 x_2^{(2)} x_3 x_1^{(q-1)}
- x_1 x_2^{(2)} x_3
- x_1 x_3^{(q)} x_2^{(2)} x_1^{(q-1)}
+ x_1 x_3^{(q)} x_2^{(2)}, \\[6pt]
&x_1 x_2^{(q)} x_1^{(q-1)}
- x_1 x_2 x_1^{(q-1)}
- x_1 x_2^{(q)}
+ x_1 x_2.
\end{align*}

Let $I_{L_2}$ be the $T$-ideal generated by $\{(x_1x_2)(x_3x_4), s_3\}$. Additionally, if $|K| = q$, we include in $I_{L_2}$ the four polynomials above as generators of this $T$-ideal.

As a consequence of the previous identities,
\(
x_1^2 x_2^{(q)} - x_1^2 x_2\) and \(x_1 x_2 x_3^{(q)} - x_1 x_2 x_3
\)
are also in $I_{L_2}$. Moreover, the polynomial
\(
x_1 x_2 x_3 x_4 - x_1 x_2 x_4 x_3
\)
is a consequence of the metabelian identity. In particular,
\(
x_1^2 x_2 x_3 \equiv_{I_{L_2}} x_1^2 x_3 x_2.
\)
Furthermore, the polynomial
\begin{equation}\label{tp}
    x_3 x_2 x_1 x_4
    - x_3 x_1 x_2 x_4
    + x_4 x_1 x_2 x_3
    - x_4 x_2 x_3 x_1
\end{equation}
is congruent to $s_3(x_1x_2, x_3, x_4)$ modulo $\langle (x_3x_4)(x_1x_2) \rangle_T$. Consequently, we also obtain the following identities in $I_{L_2}$:
\begin{equation}\label{quad}
    x_1^{(2)} x_3 x_2
    - x_1 x_3 x_1 x_2
    + x_2 x_3 x_1^{(2)}
    - x_2 x_1^{(2)} x_3
\end{equation}
and
\begin{equation}\label{quad2}
    x_1^{(2)} x_2^{(2)}
    - x_1 x_2 x_1 x_2
    + x_2^{(2)} x_1^{(2)}
    - x_2 x_1^{(2)} x_2.
\end{equation}

We have all the ingredients to prove the following result.

\begin{theorem}
Let $K$ be a finite field with cardinality $q$. Then, the polynomials
\(
(x_1 x_2)(x_3 x_4),
\)
\(s_3(x_1,x_2,x_3)= x_1x_2x_3 - x_1x_3x_2 - x_2x_1x_3 + x_2x_3x_1 + x_3x_1x_2 - x_3x_2x_1,\)
\(
x^{q+1}-x^2,
\)
and
\begin{align}
&x_1 x_2^{(2)} x_3
- x_3 x_2^{(2)} x_1 x_3^{(q-1)}
+ x_3^{(q)} x_2^{(2)} x_1
- x_1 x_3^{(q)} x_2^{(2)},\label{eq:1} \\[6pt]
&x_1 x_2^{(2)} x_3 x_1^{(q-1)}
- x_1 x_2^{(2)} x_3
- x_1 x_3^{(q)} x_2^{(2)} x_1^{(q-1)}
+ x_1 x_3^{(q)} x_2^{(2)}, \label{eq:2}\\[6pt]
&x_1 x_2^{(q)} x_1^{(q-1)}
- x_1 x_2 x_1^{(q-1)}
- x_1 x_2^{(q)}
+ x_1 x_2.\label{q2RR2}
\end{align}
form a basis of the $T$-ideal of identities of $L_2$.
\end{theorem}
\begin{proof}
First of all, let $I_{L_2}$ denote the $T$-ideal generated by all the polynomials given in our hypothesis. Let $f$ be a regular identity of $L_2$. If $f = f(x)$ depends on a single variable, then modulo $I_{L_2}$, we can write
$$
f(x) = \sum_{i=1}^{q} \lambda_i x^i.
$$
Under the evaluation $x \mapsto e_1 + \lambda e_3$, we obtain
$$
0 = f(e_1 + \lambda e_3) = \sum_{i=1}^{q} \lambda_i \alpha^{i-1}\lambda^{i-1}.
$$
In other words, any identity of $L_2$ in one variable lies in $I_{L_2}$. 

Now assume that $f = f(x,y)$ depends on two variables.  
Firstly, the identity \eqref{q2RR2} implies that $f$ contains no monomials of the form $xy$ or $yx$.  
Thus, modulo $I_{L_2}$, we can write
\[
f \;=\; 
\sum_{i,j=0}^{q-1} \lambda^0_{i,j}\, yx\, x^{(i)} y^{(j)}
\;+\;
\sum_{i,j=0}^{q-1} \lambda^1_{i,j}\, xy\, x^{(i)} y^{(j)}
\;+\;
\sum_{i=0}^{q-2} \sum_{j=1}^{q-1} \lambda^2_{i,j}\, x^2 x^{(i)} y^{(j)}
\;+\;
\sum_{i=1}^{q-1} \sum_{j=0}^{q-2} \lambda^3_{i,j}\, y^2 x^{(i)} y^{(j)},
\]
with the convention that $\lambda_{0,0}^0 = \lambda_{0,0}^1 = 0$. Notice that the evaluations $x \mapsto \lambda e_3$, $y \mapsto e_1$ and $x \mapsto e_1$, $y \mapsto \lambda e_3$ imply  $\lambda_{i,0}^0=\lambda_{0,j}^1= 0$, for all $0 < i,j  \leq q-1$.  On the other hand, under the evaluations $x \mapsto \lambda_x e_2$, $y \mapsto \lambda_y e_3$ and $x \mapsto \lambda_x e_3$, $y \mapsto \lambda_y e_2$ we can assume $\lambda_{0,j}^0=\lambda_{i,0}^1= 0$, for all $0 < i,j  \leq q-1$. Multiplying Eq.~\eqref{quad2} by $x_1^{(q-2)}$ yields
\[
    y x^{(q)} y 
    \equiv_{I_{L_2}} 
    y^{(2)} x 
    - x y^{(2)} x^{(q-1)} 
    + x^{(q)} y^{(2)}.
\]
Thus, we may also assume $\lambda_{q-1,j}^0 = 0$ for all $0<j\le q-1$.

Now, under the evaluation $x \mapsto \lambda_x e_3$, $y \mapsto e_1 + \lambda_y e_3$, for some scalars $\lambda_x, \lambda_y$ in $K$, we obtain $\lambda^3_{1,j} = 0$, for any $j=0,1,\ldots, q-2$. Similarly, $\lambda^2_{i,1} = 0$, for any $i=0,1,\ldots, q-2$. Therefore, we may assume that no monomial in $f$ is linear in either $x$ or $y$. 

By Identity~\eqref{quad2}, we may assume that $\lambda_{i,j}^3 = 0$ 
for all $i = 1, \ldots, q-1$ and $j = 0, \ldots, q-2$. 
Moreover, under the evaluation 
\[
x \mapsto \lambda_x e_3, 
\qquad 
y \mapsto e_1 + \lambda_y e_3,
\]
we see that all coefficients $\lambda^0_{i,j}$ of $f$, with $1 \leq i, j \leq q-1$, must vanish.

Next, consider the evaluation 
\[
x \mapsto e_2 + \lambda_x e_3, 
\qquad 
y \mapsto e_1 + \lambda_y e_3.
\]
In this case, any monomial beginning with $x^2$ vanishes, so the only possible nonzero terms 
must begin with $xy$. A direct computation shows that
\[
(e_2+\lambda_x e_3)(e_1+\lambda_y e_3)(e_2+\lambda_x e_3)^{(r)}(e_1+\lambda_y e_3)^{(s)}
= (-1)^{r+s+1}\lambda_x^r \lambda_y^{\,s+1} e_2,
\]
for all $2 \leq r,s \leq q-1$. Hence, all coefficients $\lambda_{i,j}^1$ vanish.

Finally, from the evaluation 
\[
x \mapsto e_1 + \lambda_x e_3, 
\qquad 
y \mapsto \lambda_y e_3,
\]
for scalars $\lambda_x, \lambda_y \in K$, we deduce that the coefficients corresponding to monomials of $f$ beginning with $x^2$ must vanish. Therefore, we conclude that $f(x,y) \in I_{L_2}$.

Now let $f(x_1,x_2,x_3)$ be a regular identity for $L_2$. Thus, we can write:
$$f = f_1 + f_2,$$
where $f_1$ contains only monomials of the form $x_i x_j^{(\epsilon_j)} x_k^{(\epsilon_k)} x_i^{(\epsilon_i)}$ and $f_2$ consists of the monomials $x_i^{(\epsilon_i)} x_j^{(\epsilon_j)} x_k^{(\epsilon_k)}$. We can then write:
$$f_1 = \sum_{(i,j,k,\epsilon_i,\epsilon_j,\epsilon_k)}\lambda_{(i,j,k)}^{(\epsilon_i,\epsilon_j,\epsilon_k)} x_i x_j^{(\epsilon_j)} x_k^{(\epsilon_k)} x_i^{(\epsilon_i-1)} $$

From the identity $s_3$, we may assume that $f_1$ contains no monomial where $x_1$ and $x_2$ are the first and second variable, respectively. Hence, if a monomial in $f_1$ starts with $x_1$, the second variable will be $x_3$. Furthermore, the identity $x_1 x_2^{(q)}x_1^{(q-1)} - x_1 x_2 x_1^{(q-1)} - x_1 x_2^{(q)} + x_1 x_2$ allows us to assume that no monomial in $f_1$ is linear in its first two variables. Therefore, all coefficients of the form $\lambda_{(1,2,k)}^{(\epsilon_1,\epsilon_2,\epsilon_k)}$ and $\lambda_{(i,j,k)}^{(1,1,\epsilon_k)}$ are zero.

Thus,
\begin{eqnarray*}
    f(e_1,\lambda_2 e_3, \lambda_3 e_3)&=&f_1(e_1,\lambda_2 e_3, \lambda_3 e_3)\\
&=&\sum_{(1,3,2,1,\epsilon_j,\epsilon_k)}\lambda_{(1,3,2)}^{(1,\epsilon_3,\epsilon_2)} e_1 (\lambda_3 e_3)^{(\epsilon_3)} (\lambda_2 e_3)^{(\epsilon_2)}\\
&=& \sum_{(1,3,2,1,\epsilon_j,\epsilon_k)}\lambda_{(1,3,2)}^{(1,\epsilon_3,\epsilon_2)} \alpha^{\epsilon_2 + \epsilon_3}  \lambda_3^{\epsilon_3} \lambda_2^{\epsilon_2} e_1,
\end{eqnarray*}
where $2 \leq \epsilon_3 \leq q$ and $1 \leq \epsilon_2 \leq q-1$. From our knowledge of the basis of identities for $K$, we have shown that no monomial in $f_1$ begins with $x_1$ and is linear in this variable. Identity \eqref{quad} then allows us to conclude that no monomial in $f_1$ has $x_1$ as its first variable.

We can also express $f_2$ as:
$$f_2 = \sum_{(i,j,k,\epsilon_i,\epsilon_j,\epsilon_k)} \theta_{(i,j,k)}^{(\epsilon_i,\epsilon_j,\epsilon_k)} x_i^{\epsilon_i} x_j^{\epsilon_j} x_k^{\epsilon_k},$$
where $\epsilon_i> 1$. Using the identity $(x_1 x_2) (x_3 x_4)$, we may assume that $j<k$ in each monomial. 
%Moreover, by $s_3$, we assume that coefficients of the form $\theta_{(1,3,2)}^{(\epsilon_1,\epsilon_3,\epsilon_2)}$ vanish. 
Therefore,
$$f(e_1 + \lambda_1 e_3, \lambda_2 e_3, \lambda_3 e_3) = \sum_{(1,2,3,\epsilon_1,\epsilon_3,\epsilon_2)} \alpha^{\epsilon_1+\epsilon_2+\epsilon_3-1}\theta_{(1,2,3)}^{(\epsilon_1,\epsilon_3,\epsilon_2)}\lambda_1^{\epsilon_1-1} \lambda_2^{\epsilon_2} \lambda_3^{\epsilon_3},$$
proving that the elements of the form $\theta_{(1,2,3)}^{(\epsilon_1,\epsilon_3,\epsilon_2)}$ are zero. Hence, $f$ contains no monomials starting with $x_1$.

Furthermore, by identities \eqref{quad} and \eqref{eq:1}, we can also assume that $\lambda_{(2,1,3)}^{(\epsilon_2,\epsilon_1,\epsilon_3)} = 0$.

Now, from identity \eqref{eq:2}, we may suppose that the coefficients of the form $\lambda_{(3,1,2)}^{(1,\epsilon_1,\epsilon_2)}$ are zero. Consequently, for $\lambda_1,\lambda_2,\lambda_3 \in K$, we have:
$$ 0 = f(e_2 + \lambda_1 e_3,  \lambda_2 e_3, \lambda_3 e_3) = \sum_{(3,1,2,\epsilon_3,\epsilon_1,\epsilon_2)}\lambda_{(3,1,2)}^{(\epsilon_3,\epsilon_1,\epsilon_2)} (-1)^{ \epsilon_1 + \epsilon_2 + \epsilon_3} \lambda_1^{\epsilon_1-1} \lambda_2^{\epsilon_2} \lambda_3^{\epsilon_3},$$
allowing us to conclude that every monomial of the form $x_3 x_1^{(\epsilon_1)} x_2^{(\epsilon_2)} x_3^{(\epsilon_3-1)}$ has coefficient $0$.

Therefore, we may write $f$ in such a way that every monomial has $x_2$ and $x_3$ as its first two variables (not necessarily in this order). Thus, for each 
$i \in \{1,\dots,q-1\}$ there exists a polynomial $g_i(x_2,x_3)$ such that
\[
    f(x_1,x_2,x_3) = \sum_{i=1}^{q-1} g_i(x_2,x_3)\, x_1^{(i)}.
\]
Now, for arbitrary elements $a_2,a_3 \in L_2$ and $\lambda_1 e_3 \in L_2$, 
consider the evaluation 
\[
    x_2 \mapsto a_2, \qquad x_3 \mapsto a_3, \qquad x_1 \mapsto \lambda_1 e_3.
\]
For each $i$ there exist scalars $\omega_{1,i},\omega_{2,i} \in K$ such that  
\[
    g_i(x_2,x_3) = \omega_{1,i} e_1 + \omega_{2,i} e_2.
\]
Hence,
\[
    f(\lambda_1 e_3, a_2, a_3)
    = \sum_{i=1}^{q-1} \lambda_1^{i} \left( \alpha^{i}\omega_{1,i} e_1 
        + (-1)^{i}\, \omega_{2,i} e_2 \right).
\]
From this, we conclude that $\omega_{1,i} = \omega_{2,i} = 0$ for all 
$i \in \{1,\dots,q-1\}$. In other words, each $g_i(x_2,x_3)$ is a polynomial identity for $L_2$, and therefore, by what has already been established, 
the polynomial $f$ follows from the identities already considered.

It is clear that the issue of variable powers becomes simpler in the case of multihomogeneous polynomials, and consequently over infinite fields. We will formalize this in the next result.

Now, let us analyze the case of a regular identity $f(x_1,x_2,x_3,\dots,x_n)$ with $n \geq 4$. We can write $f = f_1 + f_2$, where $f_1$ and $f_2$ are polynomials given by:
$$f_1(x_1,x_2,\dots,x_n) = \sum_{(i,j,\epsilon_i,\epsilon_j,\epsilon_1,\epsilon_2,\dots,\epsilon_n)} \lambda_{(i,j)}^{(\epsilon_i,\epsilon_j,\epsilon_1,\epsilon_2,\dots,\epsilon_n)} x_i x_j^{(\epsilon_j)} x_1^{(\epsilon_1)} \cdots x_n^{(\epsilon_n)}x_i^{(\epsilon_i-1)}$$
$$f_2(x_1,x_2,\dots,x_n) = \sum_{(i,\mu_i,\mu_1,\dots,\mu_n)}\theta_i^{(\mu_i,\mu_1,\dots,\mu_n)} x_i^{(\mu_i)} x_1^{(\mu_1)}x_2^{(\mu_2)}\cdots x_n^{(\mu_n)},$$
where $1 \leq\epsilon_i \leq q$, $1 \leq \epsilon_j \leq q$, $1 \leq \epsilon_k < q$, $2 \leq \mu_i \leq q$, $1 \leq \mu_\ell < q$ for $k,l \in \{1,\dots n\}$ with $k \neq i$, $k \neq j$ and $\ell \neq i$.

Take $m \in \{4,5,\dots,n\}$. The polynomial \eqref{tp} implies that the coefficients of the form $\lambda_{(i,m)}^{(\epsilon_i,\epsilon_m,\epsilon_1,\epsilon_2,\dots,\epsilon_n)}$ and $\lambda_{(m,j)}^{(\epsilon_m,\epsilon_j,\epsilon_1,\epsilon_2,\dots,\epsilon_n)}$ are zero when $i \neq 1$ and $j \neq 2$.

Evaluating $x_m \mapsto e_1$ and $x_k \mapsto \lambda_k e_3$ for $k \neq m$, we conclude that $\lambda_{(m,2)}^{(1,\epsilon_2,\epsilon_1,\epsilon_3,\dots,\epsilon_n)} = 0$.
With this, we can assume from identity \eqref{quad} that $\lambda_{(m,2)}^{(\epsilon_m,\epsilon_2,\epsilon_1,\epsilon_3,\dots,\epsilon_n)} = 0$ for any $\epsilon_1,\epsilon_2,\dots, \epsilon_n$.

Moreover, by \eqref{eq:2} and \eqref{q2RR2}, we may assume that $\lambda_{(1,m)}^{(1,\epsilon_m,\epsilon_2,\dots,\epsilon_n)} = 0$ for all $\epsilon_m,\epsilon_2,\dots,\epsilon_n$. Therefore, evaluating $x_m \mapsto e_2 + \lambda_m e_3$, $x_k \mapsto \lambda_k e_3$ for $k \neq m$, we conclude that
$$ 0 = \sum_{(1,m,\epsilon_1,\epsilon_m,\epsilon_1,\epsilon_2,\dots,\epsilon_n)} (-1)^{ \epsilon_1 + \epsilon_2 + \cdots + \epsilon_n} \lambda_{(1,m)}^{(\epsilon_1,\epsilon_m,\epsilon_1,\epsilon_2,\dots,\epsilon_n)} \lambda_m^{\epsilon_m-1} \lambda_1^{\epsilon_1} \lambda_2^{\epsilon_2} \lambda_3^{\epsilon_3} \cdots \widehat{\lambda_m^{\epsilon_m}} \cdots \lambda_n^{\epsilon_n}e_2,$$
concluding that $\lambda_{(1,m)}^{(\epsilon_1,\epsilon_m,\epsilon_1,\epsilon_2,\dots,\epsilon_n)} = 0$ regardless of the choice of $\epsilon_1,\epsilon_m,\epsilon_1,\epsilon_2,\dots,\epsilon_n$. Furthermore, the substitution $x_m \mapsto e_1 + \lambda_m e_3$, $x_k \mapsto \lambda_k e_3$ for $k \neq m$, proves that $\theta_m^{(\mu_m,\mu_1,\dots,\mu_n)} = 0$.

Since we took an arbitrary $m \in \{4,5,\dots,n\}$, we conclude that there exist polynomials $g_{\epsilon_4,\epsilon_5,\dots,\epsilon_n}(x_1,x_2,x_3)$ such that
$$f = \sum_{(\epsilon_4,\epsilon_5,\dots,\epsilon_n)} g_{\epsilon_4,\epsilon_5,\dots,\epsilon_n}(x_1,x_2,x_3) x_4^{(\epsilon_4)} x_5^{(\epsilon_5)} \cdots x_n^{(\epsilon_n)}.$$
Analogously to the previous case, we can prove that each $g_{\epsilon_4,\epsilon_5,\dots,\epsilon_n}(x_1,x_2,x_3)$ is an identity and, by the previous result, $g_{\epsilon_4,\epsilon_5,\dots,\epsilon_n}(x_1,x_2,x_3) \in I_{L_2}$. Therefore, we have shown that $f$ is an element of $I_{L_2}$, as desired.
\end{proof}

As mentioned in the proof of the previous theorem, we can establish the following result.

\begin{theorem}\label{theo:RR2-Kinf}
If $K$ is an infinite field, then the polynomials $(x_1 x_2)(x_3 x_4)$ and
\[
s_3(x_1,x_2,x_3)= x_1x_2x_3 - x_1x_3x_2 - x_2x_1x_3 + x_2x_3x_1 + x_3x_1x_2 - x_3x_2x_1
\]
generate the ideal of identities of $L_2$ as a $T$-ideal.
\end{theorem}

\begin{proof}
The proof follows \emph{mutatis mutandis} the arguments of the previous theorem, replacing the regular polynomial with the multihomogeneous one. Therefore, we omit the details.
\end{proof}

        \begin{corollary}\label{BARL L_2}
    Let $\mathcal B_{L_2} \subset \mathcal L \langle X \rangle$ be the set consisting of the monomials
    \begin{align*}
    & x_{i_1}^{(m_{i_1})}, \qquad  x_{i_1}x_{i_2}^{(m_{i_2})}x_{i_1}^{(m_{i_1}-1)},\\
    & x_{j_1}^{(m_{j_1})}x_{j_2}^{(m_{j_2})}\qquad {(\text{if } m_{j_2}> 1 \text{, then } j_1 < j_2)}\\
    &x_{i_1}x_{i_2}^{(m_{i_2})}x_{i_3}^{(m_{i_3})}\cdots x_{i_n}^{(m_{i_n})}x_{i_1}^{(m_{i_1}-1)},\\
&x_{j_1}^{(m_{j_1})}x_{j_2}^{(m_{j_2})}x_{j_3}^{(m_{j_3})}\cdots x_{j_n}^{(m_{j_n})} \quad (\text{if } j_2 < j_1 \text{ and } m_{j_3}> 1 \text{, then } j_1 < j_3),
\end{align*}
    with $j_2 < j_3 < \cdots j_n$, $i_3 < i_4< \cdots < i_n$ and $(i_1,i_2,i_3)$ is not in ascending order. Furthermore, let $W_1$ be the set of the three smallest integers in $\{i_1,i_2,\dotsc,i_n\}$ and:
    \begin{itemize}
        \item[$(i)$] if $i_1 \notin W_1$, then $i_2  < i_k$ for $k = 3,4,\dotsc,n$ and $m_{i_1} = 1$;
        \item[$(ii)$] if $i_2 \notin W_1$, then $i_1 < i_k$ for $k = 3,4,\dotsc,n$ and $m_{i_1} >1$.
    \end{itemize}

    If $K$ is an infinite field, then the images of $\mathcal B_{L_2}$ in the relatively free algebra $\mathcal L \langle X \rangle/T(L_2)$ form a basis for this algebra. Additionally, if $K$ is a finite field of cardinality $q$, then the images of $\mathcal{B}_{L_2}$ in the relatively free algebra $\mathcal{L} \langle X \rangle / T(L_2)$, subject to the restrictions:
$$
m_{i_1}, m_{i_2},m_{j_1} \leq q, \quad \text{and} \quad m_{j_2},m_{i_k},m_{j_k} < q \ \text{for all} \ k = 3, \ldots, n,
$$
\begin{itemize}
    \item If $m_{i_1} = 1$, then $m_{i_2} > 1$;
    \item In the monomials with $n=2$, we ask that if $m_{i_1} > 1$, then $m_{i_2} < q$;
    \item $i_2 < i_1 < i_3$ does not occurs;
    \item In the case of $W_1 = \{i_1,i_2,i_3\}$, the condition $i_1 < i_2,i_3$ implies $m_{i_1} = 1$;  $i_2 < i_3 < i_1$ implies $m_{i_1}> 1$; $i_3 < i_1,i_2$ and $m_{i_1}> 1$ implies $ m_{i_2}< q$.
\end{itemize}

Moreover, its codimension sequence is given by
    \[c_1(L_2) = 1, \quad c_2(L_2) = 2, \quad c_n(L_2) = 2n-1, \forall n \geq 3.\]
Furthermore, the exponent of $L_2$ is 1.
\end{corollary}

\medskip

\subsection{\texorpdfstring{The algebra $L_3$}{The algebra L3}}\

In the algebra \(L_3\), the multiplication on the basis \(\{e_1,e_2,e_3\}\) is given by  
\[
    e_3 e_3 = e_1, \qquad  e_3 e_2 = e_2, \qquad  e_2 e_3 = - e_2.
\]
Thus, let $\alpha_i, \beta_i, \gamma_i\in K$ and $a_i =  \alpha_ie_1+\beta_ie_2+\gamma_ie_3$, $i = 1, 2, 3$. Note that
    \[(\alpha_1e_1+\beta_1e_2+\gamma_1e_3)(\alpha_2e_1+\beta_2e_2+\gamma_2e_3) = (\gamma_1\beta_2-\gamma_2\beta_1)e_2 + \gamma_1\gamma_2e_1.\]

\begin{lemma}
The polynomials
\[
    x_1^2 x_2, \qquad
    x_1 x_2 x_3 + x_2 x_1 x_3, \qquad
    x_3 x_1 x_2 - x_3 x_2 x_1 - x_2 x_1 x_3,
\]
and the product \((x_1 x_2)(x_3 x_4)\) are identities for the algebra \(L_3\).  
Moreover, if the ground field \(K\) has characteristic \(2\), then the polynomial
\[
    x_1 x_2 + x_2 x_1
\]
is also an identity for \(L_3\).

Additionally, if \(K\) is a finite field with \(|K| = q\), then the following polynomials are identities for \(L_3\):
\[
    x_1 x_2 x_3^{(q)} - x_1 x_2 x_3,
\]
\[
    x_1 x_2 - x_2 x_1
    + 2 x_2 x_1 x_2^{(q-1)}
    + 2 x_2 x_1^{(q)}
    - 2 x_2 x_1^{(q)} x_2^{(q-1)},
\]
\[
    x_3 x_1 x_2
    - x_3 x_1^{(q)} x_2
    - x_3 x_1 x_3^{(q-1)} x_2
    + x_3 x_1^{(q)} x_3^{(q-1)} x_2.
\]
\end{lemma}
\begin{proof} The verification is routine: each polynomial evaluates to zero under all substitutions of elements in $L_3$. 
\end{proof}

Notably, the identity \(x_1^2x_2\) follows from \(x_1x_2x_3 + x_2x_1x_3\) and \(x_3x_1x_2 - x_3x_2x_1 - x_2x_1x_3\) by substituting \(x_3\) by \(x_1\). Moreover, in characteristic 2,  the identity $x_1x_2x_3  + x_2x_1x_3$  follows from $x_1x_2 + x_2x_1$ by right multiplication by \(x_3\).

Let $I_{L_3}$ the $T$-ideal generated by the identities of the previous lemmas. In general, we have the following equivalence
$$x_1x_2 \equiv_{I_{L_3}} x_2x_1 - 2x_2x_1x_2^{(q-1)} - 2x_2x_1^{(q)} + 2x_2x_1^{(q)}x_2^{(q-1)}.$$

\begin{corollary}\label{idg2RR3}
    In a finite field with cardinality $|K| = q$, every identity in two variables follows from the identities $x_1x_2x_3 + x_2x_1x_3$, $x_3x_1x_2 - x_3x_2x_1 - x_2x_1x_3$, $(x_1x_2)(x_3x_4)$, $x_1x_2x_3^{(q)} - x_1x_2x_3$ and $x_1x_2 - x_2x_1 + 2x_2x_1x_2^{(q-1)} + 2x_2x_1^{(q)} -2x_2x_1^{(q)}x_2^{(q-1)}$.
\end{corollary}
\begin{proof}
    Let $f$ be a regular identity of $L_3$, then
    \[f\equiv_{I_{L_3}} \sum_{\mbox{\tiny $\begin{array}{c}
        1\leq m_1\leq q \\
        0\leq m_2 < q 
\end{array}$}}\lambda_{(m_1,m_2)}x_2x_1^{(m_1)}x_2^{(m_2)}.\]
  Initially, note that
    \(f(e_3, e_3) = \lambda_{(1,0)}e_1,\)
    then $\lambda_{(1,0)} = 0$. Thus, considering the evaluation $\varphi: x_1\mapsto \beta_1e_3, x_2\mapsto \beta_2e_2$, we have
    \[\left(
    \sum_{\mbox{\tiny $\begin{array}{c}
     2 \leq m_1\leq q 
    \end{array}$}}\lambda_{(m_1,0)}(-1)^{m_1}\beta_1^{m_1}\beta_2 \right)e_2 = 0.\]
    So, from the canonical basis of $K\langle X\rangle/T(K)$ we have that $\lambda_{(m_1,0)} = 0$ for $1 \leq m_1 \leq q$. Then, we can evaluate $x_1 \to e_2 + \beta_1 e_3$ and $x_2 \to \beta_2 e_3$ obtaining:
    \[\left(
    \sum_{\mbox{\tiny $\begin{array}{c}
        1 \leq m_1\leq q \\
        1 \leq m_2 < q
    \end{array}$}}\lambda_{(m_1,m_2)}(-1)^{m_1+m_2-1}\beta_1^{m_1-1}\beta_2^{m_2+1} \right)e_2 = 0.\]
    By the same argument, we have that $\lambda_{(m_1,m_2)} = 0$ for $1 \leq m_1 \leq q$ and $1 \leq m_2 < q$. With this, we conclude that $f \in I_{L_3}$.
\end{proof}

Now we can state the main result for the identities of $L_3$.

\begin{theorem}
Let \(K\) be an infinite field.
\begin{itemize}
    \item[$(1)$] If \(\operatorname{char}(K)\neq 2\), then the polynomials
    \[
        x_1 x_2 x_3 + x_2 x_1 x_3,\qquad
        x_3 x_1 x_2 - x_3 x_2 x_1 - x_2 x_1 x_3,\qquad
        (x_1 x_2)(x_3 x_4)
    \]
    generate the \(T\)-ideal of polynomial identities of \(L_3\).

    \item[$(2)$] If \(\operatorname{char}(K)=2\), then the polynomials
    \[
        x_1 x_2 + x_2 x_1,\qquad
        (x_1 x_2)(x_3 x_4)
    \]
    generate the \(T\)-ideal of identities of \(L_3\).
\end{itemize}
Furthermore, if \(K\) is a finite field of cardinality \(|K| = q < \infty\), then one must additionally include the identities
\[
    x_1 x_2 - x_2 x_1
    + 2 x_2 x_1 x_2^{(q-1)}
    + 2 x_2 x_1^{(q)}
    - 2 x_2 x_1^{(q)} x_2^{(q-1)},
\]
\[
    x_1 x_2 x_3^{(q)} - x_1 x_2 x_3,
\]
\[
    x_3 x_1 x_2
    - x_3 x_1^{(q)} x_2
    - x_3 x_1 x_3^{(q-1)} x_2
    + x_3 x_1^{(q)} x_3^{(q-1)} x_2,
\]
to the corresponding generating sets above.
\end{theorem}

\begin{proof}
First, denote by \(I_{L_3^{(1)}}\) and \(I_{L_3^{(2)}}\) the \(T\)-ideals generated by the sets of identities listed in items \((1)\) and \((2)\), respectively.  
For the sake of notation, in what follows we write \(I = I_{L_3^{(1)}}\) or \(I = I_{L_3^{(2)}}\), according to the characteristic of the field.

Over a field of characteristic \(2\), the polynomial
\[
    x_3 x_1 x_2 - x_3 x_2 x_1 - x_2 x_1 x_3
\]
is a consequence of the Leibniz identity together with the relation \(x_1 x_2 + x_2 x_1\).

Now, let \(f = f(x_1,x_2,\dots,x_n)\) be a multihomogeneous identity of \(L_3\).  
Since any identity of degree \(1\) or \(2\) already lies in \(I\), we may assume that the length of $f > 2$.  
Modulo \(I\), we may write \(f\) in the form
\[
    f \equiv_I 
    \sum_{j=1}^{n-1}
        \lambda_j\,
        x_n x_j^{(m_j)} x_1^{(m_1)} \cdots 
        \widehat{x_j^{(m_j)}} \cdots 
        x_n^{(m_n-1)},
\]
where the notation \(\widehat{x_j^{(m_j)}}\) indicates omission of that factor and the coefficients \(\lambda_j \in K\).

Consider now the evaluation $x_j\mapsto e_2 + e_3$ and $x_k\mapsto e_3$, we obtain \[\lambda_j(-1)^{\sum_{j=1}^{n}m_j}e_2 = 0 \Rightarrow \lambda_j = 0, \forall j\in \{1, \dotsc, n-1\}.\] 
Therefore, $f$ belongs to \(I\). This shows that \(I\) contains every multihomogeneous identity of \(L_3\). This completes the argument for the case when the ground field is infinite.

To conclude, suppose that $K$ is a finite field of cardinality $|K| = q$.  
Let $f = f(x_1,\dotsc,x_n)$ be a regular identity of $L_3$.  
By Corollary~\ref{idg2RR3}, we may assume that $n \ge 3$.  
Moreover, we can write
\[
    f \equiv_I 
    \sum_{j}
    \lambda_j \,
    x_n x_{r_j}^{(m_{r_j j})}
    x_1^{(m_{1j})} \cdots 
    \widehat{x_{r_j}^{(m_{r_j j})}} \cdots 
    x_n^{(m_{nj})}.
\]
The identity $x_1 x_2 x_3^{(q)} - x_1 x_2 x_3$ allows us to impose the restrictions
\[
    1 \le m_{r_j j} \le q,\qquad
    0 \le m_{n j} < q,\qquad
    1 \le m_{i j} < q \;\; \text{for } i \ne r_j,n.
\]
  Let $J = J_0\cup J_{q-1}$, where $J_0=\{ j \mid m_{nj} = 0 \}$ and $J_{q-1}=\{ j \mid m_{nj} = q-1\}$ and fix $j$. Considering the evaluation $x_{r_j}\mapsto e_2 +\beta_{r_j} e_3$ and $x_{r_j}\mapsto e_2 +\beta_k e_3$ if $k\neq r_j$, we can assume $\lambda_j = 0$ for all $j\notin J$ and $J_0$ and $J_{q-1}$ both no empty. Notice that if a term of the form $x_n x_{r_j}^{(q)} x_r$ appears for some 
$j \in J_0$, then, using the equivalence
\[
    x_n x_{r_j}^{(q)} x_r 
    \equiv_I 
    x_n x_{r_j} x_r 
    - x_n x_{r_j} x_n^{(q-1)} x_r 
    + x_n x_{r_j}^{(q)} x_n^{(q-1)} x_r,
\]
we may assume that such a monomial is linear in $x_{r_j}$.  Now consider the evaluation  
\[
    x_n \mapsto \beta_n e_2, \qquad 
    x_k \mapsto \beta_k e_3 \quad (1 \le k \le n-1),
\]
where $\beta_1,\dots,\beta_n \in K$ are arbitrary.  
Under this substitution we obtain
\[
    \sum_{j \in J_0} 
        \lambda_j 
        (-1)^{\sum_{i=1}^{n-1} m_{ij}}
        \beta_1^{m_{1j}}
        \cdots 
        \beta_{n-1}^{m_{(n-1)j}}
        \beta_n 
    = 0.
\]
Since the $\beta_i$ are arbitrary, it follows that 
$\lambda_j = 0$ for all $j \in J_0$.  
Consequently, $\lambda_j = 0$ for all $j \in J_{q-1}$ as well, 
and the result follows.
\end{proof}

\begin{corollary}\label{BARL L_3}
    Let $\mathcal{B}_{L_3} \subset \mathcal{L}\langle X\rangle$ be the set formed by the monomials
    \[
        x_i,\qquad 
        x_j x_k,\qquad
        x_{j_1} x_{j_2}^{(m_2)} \cdots x_{j_n}^{(m_n)}
        \quad (n \ge 2,\ j_1 > j_2,\ j_3 < j_4 < \dots < j_n \le j_1,\ 
        j_2 \neq j_k\ \forall k>2).
    \]
    Suppose that $K$ is an infinite field.  
    If $\operatorname{char} K \ne 2$, then the images of the monomials of 
    $\mathcal{B}_{L_3}$ form a basis of the relatively free algebra 
    $\mathcal{L}\langle X\rangle / T(L_3)$. When $|K|=q<\infty$, the same set forms a basis subject to the additional restrictions:
    \begin{itemize}
        \item $1 \le m_2 \le q$, $0 \le m_n < q$, and 
              $1 \le m_k < q$ for all $k \ne 2,n$;

        \item if $m_2 = q$, then $m_n > 0$;

        \item if $m_2 = 1$ and $m_n = 0$, then $j_2 < j_3$.
    \end{itemize}

    Furthermore,
    \[
        c_1(L_3) = 1,\qquad 
        c_2(L_3) = 2,\qquad 
        c_n(L_3) = n-1 \quad (n>2).
    \]

    If $\operatorname{char} K = 2$, then one must additionally impose the condition 
    $j \ge k$ for the quadratic monomials $x_j x_k$, and in this case $c_2(L_3)=1$.

    In all cases, the exponent of $L_3$ is equal to $1$.
\end{corollary}

\medskip

\subsection{\texorpdfstring{The algebra $L_9$}{The algebra L9}}\

In the algebra \(L_9\), the multiplication on the basis \(\{e_1, e_2, e_3\}\) is given by  
\[
    e_1 e_3 = e_2, \qquad  e_2 e_3 = \alpha e_1.
\]

Now, for arbitrary scalars \(\alpha_i, \beta_i, \gamma_i \in K\) with \(i = 1, \dots, n\), it follows by a straightforward computation that, if \(n\) is even, then
\[
(\alpha_1e_1+\beta_1e_2+\gamma_1e_3)\cdots(\alpha_ne_1+\beta_ne_2+\gamma_ne_3)
= \alpha^{\frac{n}{2}-1}\gamma_2\cdots\gamma_n\,(\alpha_1 e_2 + \alpha \beta_1 e_1),
\]
while if \(n\) is odd, we obtain
\[
(\alpha_1e_1+\beta_1e_2+\gamma_1e_3)\cdots(\alpha_ne_1+\beta_ne_2+\gamma_ne_3)
= \alpha^{\frac{n-1}{2}}\gamma_2\cdots\gamma_n\,(\alpha_1 e_1 + \beta_1 e_2).
\]

Thus, in particular, the polynomial \(x_1(x_2x_3)\) is an identity for \(L_{9}\).

Furthermore, suppose that \(|K| = q < \infty\). Then we have two cases: If \(q\) is even, then  
\(
x_1 x_2^{(2q-1)} - x_1 x_2\) and \(x_1 x_2^{(q)} x_3 - x_1 x_2 x_3^{(q)}\) are in $T(L_{9})$. If \(q\) is odd, then  
\(
x_1 x_2^{(q)} - \alpha^{\frac{q-1}{2}}\, x_1 x_2\) lies in $T(L_{9})$.

\begin{theorem}
    If $K$ is an infinite field, then the polynomial $x_1(x_2x_3)$ generates the $T$-ideal of the identities of $L_{9}$. 
    If $K$ is a finite field with $|K| = q$, then:
    \begin{itemize}
        \item[$(i)$] If $q$ is odd, the polynomials 
        \[
            x_1(x_2x_3), \qquad x_1x_2^{(q)} - \alpha^{\frac{q-1}{2}} x_1x_2
        \]
        generate the $T$-ideal of the identities of $L_{9}$;
        
        \item[$(ii)$] If $q$ is even, the polynomials
        \[
            x_1(x_2x_3), \qquad x_1 x_2^{(2q-1)} - x_1 x_2,\qquad 
            x_1 x_2^{(q)} x_3 - x_1 x_2 x_3^{(q)}
        \]
        generate the $T$-ideal of the identities of $L_{9}$.
    \end{itemize}
\end{theorem}
\begin{proof}
Let $f = f(x_1,\dots,x_n)$ be a multihomogeneous identity of $L_{9}$.  
Modulo the $T$-ideal generated by $x_1(x_2x_3)$, we may write
\[
    f \equiv
    \sum_{j=1}^{n} 
        \lambda_j\,
        x_j^{(m_j)}
        x_1^{(m_1)}
        \cdots
        \widehat{x_j^{(m_j)}}
        \cdots
        x_n^{(m_n)}.
\]
Now apply the evaluation 
\[
    x_j \mapsto e_1 + e_3,\qquad 
    x_k \mapsto e_3 \ \text{for } k \ne j.
\]
Since $f$ is an identity, the evaluation must vanish, forcing $\lambda_j = 0$ for all $j$. This proves that every identity follows from $x_1(x_2x_3)$ when $K$ is infinite.

Suppose that $K$ is a finite field with cardinality $|K| = q$. Let $f$ be a regular identity of $L_9$, not necessarily multihomogeneous. We begin with the case where $q$ is odd. 

Let $I_{L_9}$ be the $T$-ideal generated by the identities in the hypothesis. We may then express
\begin{equation}\label{idrRR13}
f \equiv_{I_{L_9}} \sum_{j} \lambda_j \, x_{r_j}^{(m_{r_j j})} x_1^{(m_{1 j})} \cdots \widehat{x_{r_j}^{(m_{r_j j})}} \cdots x_n^{(m_{n j})}.
\end{equation}
Now, for each $k \in \{1, \dotsc, n\}$, define
\(
J_k = \{ j \mid r_j = k \}.
\)
Under the evaluation
\[
x_k \mapsto e_1 + \gamma_k e_3, \qquad x_i \mapsto \gamma_i e_3 \ \text{for } i \neq k,
\]
with arbitrary scalars $\gamma_i \in K$, the image of each term in \eqref{idrRR13} is
\[
\sum_{j \in J_k} \lambda_j \, \gamma_k^{m_{kj}-1} \gamma_1^{m_{1j}} \cdots \widehat{\gamma_k^{m_{r_j j}}} \cdots \gamma_n^{m_{nj}} \, u(j) = 0,
\]
where
\[
m_j = m_{1j} + m_{2j} + \dotsb + m_{nj}, \qquad
u(j) =
\begin{cases}
\alpha^{\frac{m_j-1}{2}} e_1, & \text{if $m_j$ is odd},\\[1ex]
\alpha^{\frac{m_j}{2}-1} e_2, & \text{if $m_j$ is even}.
\end{cases}
\]

Since the $\gamma_i$ are arbitrary and from the bounds on the powers of the variables imposed by the finite field, we conclude that $\lambda_j = 0$ for all $j \in J_k$. Repeating this argument for all $k = 1, \dotsc, n$, it follows that $f \equiv 0$ modulo the $T$-ideal generated by the prescribed identities. This completes the proof in the case of a finite field with odd cardinality.

  For $q$ even, denote by $I_{L_{9}}$ the $T$-ideal generated by the identities corresponding to this case. In this situation, $f$ can be written in the form given by \eqref{idrRR13}. Fix $k \in \{1,\ldots,n\}$ and define
\[
J_k = \{\, j \mid r_j = k,\ m_{kj} = 1 \,\}.
\]
Inside $J_k$, consider the subsets
\[
J_k^{1} = \{\, j \in J_k \mid m_j \text{ is odd}\},
\qquad
J_k^{2} = \{\, j \in J_k \mid m_j \text{ is even}\}.
\]
Clearly, $J_k = J_k^{1} \cup J_k^{2}$. Moreover, evaluating
\[
x_k \mapsto e_1,
\qquad 
x_i \mapsto \gamma_i e_3 \ \text{for } i \ne k,
\]
we obtain
\[
\sum_{j \in J_k}
\lambda_j \,\gamma_1^{m_{1j}}\cdots
\widehat{\gamma_k^{\,m_{kj}}}\cdots 
\gamma_n^{m_{nj}}\, u(j) = 0.
\]
Fix $t \in \{1,\ldots,n\}\setminus\{k\}$. By the identities 
\[
x_1 x_2^{(2q-1)} - x_1 x_2,
\qquad
x_1 x_2^{(q)} x_3 - x_1 x_2 x_3^{(q)},
\]
we may assume that $m_{tj} < 2q$ and $m_{ij} < q$ whenever $i \neq k,t$.

The previous equality, together with the definition of $u(j)$, yields
\[
\sum_{j \in J_k^{1}}
\lambda_{j}\,\alpha^{\frac{m_j-1}{2}}\,
\gamma_1^{m_{1j}}\cdots 
\widehat{\gamma_k^{\,m_{kj}}}\cdots 
\gamma_n^{m_{nj}}
\;=\; 0
\;=\;
\sum_{j \in J_k^{2}}
\lambda_{j}\,\alpha^{\frac{m_j}{2}-1}\,
\gamma_1^{m_{1j}}\cdots 
\widehat{\gamma_k^{\,m_{kj}}}\cdots 
\gamma_n^{m_{nj}}.
\]

Note that, in order for two powers in $f$ to produce the same polynomial in the commutative free algebra $K[\gamma_i \mid i > 0]$ after evaluation, their exponents must satisfy that $q+s$ and $s+1$ are congruent modulo $q-1$ and have opposite parities, with $0 \le s \le q-2$. Thus, if $j_1, j_2 \in J_k^{1} \cup J_k^{2}$ are such that
\[
m_{t j_1} \equiv m_{t j_2} \pmod{q-1}
\qquad\text{and}\qquad
m_{i j_1} = m_{i j_2} \text{ for all } i \ne t,
\]
then we may assume $j_1 \in J_k^{1}$ and $j_2 \in J_k^{2}$. It follows that
\[
\lambda_j = 0 \qquad \text{for all } j \in J_k^{1} \cup J_k^{2},
\]
that is, the evaluation forces each corresponding scalar coefficient to vanish.

Finally, in \eqref{idrRR13} the exponents satisfy the restrictions $1 < m_{r_j j} < 2q$ and $m_{ij} < q$ for all $i \ne r_j$. Evaluating
\[
x_k \mapsto e_1 + \gamma_k e_3, 
\qquad
x_i \mapsto \gamma_i e_3 \quad (i \ne k),
\]
we obtain the same type of conclusion as in the previous case, thereby completing the proof.
\end{proof}

We conclude this subsection with the following result.

\begin{corollary}
    Let $\mathcal B_{L_{9}} \subset \mathcal L \langle X \rangle$ be the set consisting of the monomials
    \[
        x_{j_1}^{(m_1)} x_{j_2}^{(m_2)} \cdots x_{j_n}^{(m_n)}
        \qquad
        (j_2 < j_3 < \dotsb < j_n,\; j_1 \neq j_k\ \text{for all } k \ge 2).
    \]
    If $K$ is an infinite field, then the images of $\mathcal B_{L_{9}}$ in the relatively free algebra
    \[
        \mathcal L \langle X \rangle / T(L_{9})
    \]
    form a basis of this algebra.

    If $K$ is finite with $|K| = q$ odd, we may impose the restrictions
    \[
        1 \le m_1 \le q,\qquad 1 \le m_i < q\ \text{for all } i \ge 2.
    \]
    In the case where $q$ is even, we may assume
    \[
        1 \le m_1 < 2q,\qquad
        1 \le m_2 < 2q-1,\qquad
        1 \le m_i < q\ \text{for all } i > 2,
    \]
    and, moreover, when $m_1 > 1$ we additionally require $m_2 < q$.

    In all cases, the codimension sequence of $L_{9}$ is
    \(
        c_n(L_{9}) = n.
    \)
    Furthermore, the exponent of $L_9$ is equal to $1$.
\end{corollary}

\medskip

\subsection{\texorpdfstring{The algebra $L_{10}$}{The algebra L10}}\

In the algebra $L_{10}$, the multiplication on the basis $\{e_1,e_2,e_3\}$ is given by  
\[
    e_1 e_3 = e_2, \qquad  
    e_2 e_3 = \alpha e_1 + e_2 .
\]
Thus, the product for arbitrary elements  
\(  x = \alpha_1 e_1 + \beta_1 e_2 + \gamma_1 e_3\) and \( y = \alpha_2 e_1 + \beta_2 e_2 + \gamma_2 e_3\),
in $L_{10}$ is
\[
    xy
    = \gamma_2\bigl(\alpha\beta_1\, e_1 + (\alpha_1 + \beta_1)\, e_2 \bigr).
\]
Consequently, the polynomial
\begin{equation}\label{RR7fin1}
    x_1(x_2 x_3)
\end{equation}
is an identity of $L_{10}$.

Moreover, as an immediate consequence of the Leibniz identity~\eqref{ideLeibniz}, the polynomial  
\(
    x_1 x_2 x_3 - x_1 x_3 x_2
\)
also lies in the $T$-ideal $T(L_{10})$.

Now consider the sequence \(\{\theta_n\}_{n\ge -1}\) defined by
\[
    \theta_{-1} = 0,\qquad
    \theta_0 = 1,
\]
and, for any positive integer \(n \ge 0\),
\[
    \theta_{n+1} = \alpha\, \theta_{n-1} + \theta_{n}.
\]
It is easy to prove, by induction, that for every integer \(m>1\),
\begin{equation}\label{mulittableRR7}
    (\alpha_1 e_1 + \beta_1 e_2 + \gamma_1 e_3)\, e_3^{\,(m)}
    = \alpha\bigl(\alpha_1 \theta_{m-2} + \beta_1 \theta_{m-1}\bigr)e_1
      + \bigl(\alpha_1 \theta_{m-1} + \beta_1 \theta_{m}\bigr)e_2 .
\end{equation}

From the multiplication table we observe that, while the first factor may be arbitrary, 
in the second factor only the component along \(e_3\) affects the product.  
Hence right multiplication by \(e_3\) induces a linear transformation on the subspace
\(\mathrm{span}\{ e_1, e_2\}\), which is represented, with respect to the ordered basis \((e_1,e_2)\), by the matrix
\[
    A=\begin{pmatrix} 0 & \alpha \\[4pt] 1 & 1 \end{pmatrix}.
\]

Moreover, by direct computation one verifies that, for all $m\ge1$,
\[
    A^m
    =
    \begin{pmatrix}
\alpha\,\theta_{m-2} & \alpha\,\theta_{m-1}\\[4pt]
\theta_{m-1 } & \theta_{m}
\end{pmatrix}.
\]

We are interested in studying polynomial identities on the algebra $L_{10}$. The next result allows us to restrict our attention to identities that begin with the same variable.

\begin{lemma}\label{startsame}
    Any identity of \( L_{10} \) is a linear combination of polynomial identities that all start with the same variable.
\end{lemma}
\begin{proof}
    Let \( f(x_1, \dots, x_s) \) be a polynomial identity of \(L_{10}\). 
    We decompose \(f\) as
    \[
        f(x_1, \dots, x_s)
        = g(x_1, \dots, x_s) + h(x_1, \dots, x_s),
    \]
    where \(g\) is the sum of all monomials of \(f\) whose leftmost variable is \(x_1\), and
    \(h\) consists of the remaining monomials.

    Consider the evaluation \(\varphi\colon \mathcal L\langle X\rangle \to L_{10}\) defined by sending
    \(x_1\) to an arbitrary element of \(L_{10}\), and
    \[
        x_i \mapsto \lambda_i e_3, \qquad i \ge 2,
    \]
    where each \(\lambda_i\in K\).  
    Since any product in \(L_{10}\) in which the leftmost factor is a scalar multiple of \(e_3\) equals zero, every monomial in \(h\) vanishes under~\(\varphi\). Thus
    \[
        0=\varphi\bigl(f(x_1,\dots,x_s)\bigr)
        =\varphi\bigl(g(x_1,\dots,x_s)\bigr).
    \]
    Hence \(g\) is itself a polynomial identity of \(L_{10}\).

    Because the choice of \(x_1\) was arbitrary, the same argument applies to any variable:
    decomposing \(f\) according to its leading variable shows that every identity of \(L_{10}\)
    is a linear combination of identities whose monomials all start with the same variable.
\end{proof}

Before continuing, we assume that the ground field $K$ is finite, with cardinality $|K|=q$.  
Under this assumption, the identity
\begin{equation}\label{RR7fin2}
    x_1 x_2^{\,(q)} x_3 - x_1 x_2 x_3^{\,(q)}
\end{equation}
holds in $T(L_{10})$, because of the variables $x_2$ and $x_3$ can only be evaluated as scalar multiples of $e_3$.  
Moreover, the product defined in \eqref{mulittableRR7} does not depend on the specific scalars in $x_2$ or $x_3$, but only on the number of occurrences of the vector $e_3$ on the right-hand side of $x_1$.  
Finally, since $\gamma^q = \gamma$ for all $\gamma \in K$, the identity follows immediately.

From the previous identity, it follows that every monomial that is not linear in the variable appearing in the first position is congruent, modulo the $T$-ideal generated by \eqref{RR7fin1} and \eqref{RR7fin2}, to a monomial of the form
\begin{equation}\label{m1}
    x_1^{(m_1)} x_2^{(m_2)} \cdots x_n^{(m_n)}, \quad \text{with } m_i < q \ \text{for } i \ge 2.
\end{equation}
On the other hand, if the variable in the first position is linear, the monomial takes the form
\begin{equation}\label{m2}
    x_1 x_2^{(m_2)} \cdots x_n^{(m_n)}, \quad \text{with } m_i < q \ \text{for } i > 2.
\end{equation}

The next step is to further restrict the degree of the first variable in monomials of type \eqref{m1}, and, for monomials of type \eqref{m2}, to bound the degree of the second variable.

We now analyze the order of the matrix \(A\) over the ground field under the above assumptions.

Since \(\alpha \neq 0\), the iterative process depends only on the initial term and the number of iterations. As \(\mathrm{GL}(2,K)\) is finite, the matrix \(A\in \mathrm{GL}(2,K)\) has finite order, which we denote by \(\operatorname{ord}(A)=\mathfrak{o}\). Consequently,
\[
A^{\mathfrak{o}}e_i=e_i,
\]
for \(i=1,2\).

The Jordan decomposition of \( A \) is of the form
\[
    \begin{pmatrix} r_1 & 1 \\[4pt] 0 & r_1 \end{pmatrix}
    \qquad\text{or}\qquad
    \begin{pmatrix} r_1 & 0 \\[4pt] 0 & r_2 \end{pmatrix},
\]
where \( r_1 \) and \( r_2 \) are the roots of the polynomial $ p(t) = t^{2} - t - \alpha$. Thus, the eigenvalues of \( A \) depend on the choice of the parameter \( \alpha \). Moreover, by Galois theory, if \( r_1 \) is a root of \( p(t) \), then its conjugate under the Frobenius automorphism is
$r_2 = r_1^{\,q}$.

In summary:
\begin{itemize}
    \item If \( r_1 \) and \( r_2 \) are two distinct roots in \( K \), then $\operatorname{ord}(A)
        = \operatorname{lcm}\bigl( \operatorname{ord}(r_1), \operatorname{ord}(r_2) \bigr)$.
    \item If \( r_1 = r_2 = r \), then
\[
    A \sim 
    \begin{pmatrix}
        r & 1 \\[4pt]
        0 & r
    \end{pmatrix}=\tilde{A},
\]
since it is not a scalar matrix. Furthermore, by Girard’s relations, the condition that \( r \) is a double
root of \( p(t)=t^{2}-t-\alpha \) implies \( 2r = 1 \), and therefore
\[
    r = \frac{1}{2},
    \qquad
    \alpha = -\frac{1}{4}.
\]
This situation can occur only if \( \operatorname{char}(K)=p\neq 2 \),
since \( 2^{-1} \) must exist in \( K \).

For such a Jordan block we have
\begin{equation}\label{eqr=r}
    \tilde{A}^{m}
    = r^{m} I + m r^{m-1} N,
    \qquad N^{2}=0.
\end{equation}
Thus \( \tilde{A}^{m}=I \) if and only if
\[
    r^{m} = 1
    \quad\text{and}\quad
    m \equiv 0 \pmod{p}.
\]
Consequently,
\(
    \operatorname{ord}(A)
        = \operatorname{lcm}\bigl( p,\ \operatorname{ord}(r) \bigr)
        = \operatorname{lcm}\bigl( p,\ \operatorname{ord}(2) \bigr),
\)
since \( r = 1/2 \) and \( \operatorname{ord}(1/2)=\operatorname{ord}(2) \).

Finally, since \( \operatorname{ord}(2) < p \) in any finite field of characteristic \(p\) (Fermat's little theorem),
we obtain the estimate
\[
     \operatorname{ord}(A)
     = \operatorname{lcm}\bigl(p,\operatorname{ord}(2)\bigr)
     < p^{2}.
\]

    \item If \( r_1, r_2 \notin K \) but lie in the quadratic extension \( F \), then
    \[
        \operatorname{ord}(A)
        = \operatorname{lcm}\bigl( \operatorname{ord}(r_1), \operatorname{ord}(r_1^{q}) \bigr)
        = \operatorname{ord}(r_1),
    \]
    because \( r_1 \) and \( r_1^{q} \) have the same multiplicative order.
\end{itemize}

Due to the previous remarks, we are now ready to prove the following result.

\begin{lemma}
Let $K$ be a finite field with cardinality $|K| = q$, and let $\alpha \in K$ be the scalar arising from the multiplication table of $L_{10}$ such that the two roots $r_1,r_2$ of 
\[
    p(t)=t^{2}-t-\alpha
\]
lie in $K$. Under these hypotheses, we have:

\begin{enumerate}
    \item If $r_1$ and $r_2$ are distinct elements of $K$, then the polynomial
\begin{equation}\label{RR7cota}
    x_1 x_2^{(q)} - x_1 x_2
\end{equation}
is an identity of $L_{10}$. Moreover, the identity \eqref{RR7fin2} follows as a consequence of this one.

    \item If $\operatorname{char}(K)=p\neq 2$ and $r_1=r_2=r$, then the polynomial
\begin{equation}\label{RR7fin4} 
x_1 x_2^{(2q-1)} - 2 x_1 x_2^{(q)} + x_1 x_2
\end{equation}    
   is an identity of $L_{10}$.
\end{enumerate}
\end{lemma}
\begin{proof}
Denote by \(\widetilde{A}\) the Jordan matrix associated to \(A\). Then there exists a matrix \(P \in \mathrm{GL}(2,K)\) such that \(
\widetilde{A} = P^{-1} A P.
\)

\noindent
\textbf{(1) Distinct roots.}
If \(r_1\neq r_2\) lie in \(K\), then the matrix \(A\) is diagonalizable over \(K\). Hence,
\[
    A \sim 
    \widetilde{A}
    = \begin{pmatrix}
        r_1 & 0 \\[4pt]
        0 & r_2
      \end{pmatrix}.
\]
Since the field \(K\) has \(q\) elements, every scalar \(\lambda\in K\) satisfies
\(\lambda^{q} = \lambda\).  
Applying this fact to the eigenvalues of \(A\), we obtain
\[
    \widetilde{A}^{\,q}=\widetilde{A}.
\]
Because \(A\) is similar to \(\widetilde{A}\), it follows that
\(A^{q}=A\). Therefore, in the algebra \(L_{10}\) the identity
\(
    x_1 x_2^{(q)} = x_1 x_2
\)
holds. The second assertion follows immediately.

\smallskip

\noindent
\textbf{(2) Double root case.}
If $r_1 = r_2 = r$, then the Jordan form of $A$ is
\[
A \sim 
\widetilde{A}
=
\begin{pmatrix}
    r & 1 \\[4pt]
    0 & r
\end{pmatrix}.
\]
From the previous discussion we require that $\operatorname{char}(K)=p \ne 2$.

Using the standard formula given in \eqref{eqr=r}, the fact that $r^{q}=r$ for all $r \in K$, and that $q \equiv 0 \pmod{p}$, we obtain
\[
\widetilde{A} = r I + N, 
\qquad 
\widetilde{A}^{q} = r I,
\qquad 
\widetilde{A}^{2q-1} = r I - N.
\]
Thus,
\[
 A^{2q-1} - 2 A^{\,q} + A = 0.
 \]
Therefore, for all $\lambda \in K$,
\[
e_1 (\lambda e_3)^{(2q-1)}
- 2\, e_1 (\lambda e_3)^{(q)}
+ e_1 (\lambda e_3)
=0,
\]
and the same computation holds with $e_2$ in place of $e_1$.
This proves that the polynomial
\[
    x_1 x_2^{(2q-1)} - 2 x_1 x_2^{(q)} + x_1 x_2
\]
This completes the proof of both cases.
\end{proof}

We now turn to the third case, namely when the two roots $r_1$ and $r_2$ of 
\(
    p(t)=t^2 - t - \alpha
\)
are distinct but do not lie in the base field $K$.  
Equivalently, $r_1 \notin K$ (and hence $r_2 = r_1^{\,q}$), so the splitting field of $p(t)$ is the quadratic extension $F$ of $K$.

The next result provides sharper bounds for the behavior of the first two variables in this situation, where the Frobenius conjugates $r_1$ and $r_2$ lie in $F \setminus K$.

\begin{lemma}\label{limitaçãoRR7finite}
Let $K$ be a finite field with cardinality $|K| = q$, and let $\alpha$ be a scalar from the multiplication table of $L_{10}$ such that the roots $r_1$ and $r_2$ of the polynomial
\(
p(t)=t^2 - t - \alpha
\)
do not lie in $K$. Then
\[
    r_1^{\,q-1} + r_2^{\,q-1} \in K,
\]
and the polynomial
\begin{equation}\label{RR7fin3}
    x_1 x_2^{(2q-1)}
    \;-\;
    \bigl( r_1^{\,q-1} + r_2^{\,q-1} \bigr) x_1 x_2^{(q)}
    \;+\;
    x_1 x_2
\end{equation}
% and
% \begin{equation}\label{RR7fin6}
%     x_1^{(2q)} x_2
%     \;-\;
%     \bigl( r_1^{\,q-1} + r_2^{\,q-1} \bigr) x_1^{(q+1)} x_2
%     \;+\;
%     x_1^{(2)} x_2
% \end{equation}
is an identity of the algebra \(L_{10}\).
\end{lemma}
\begin{proof}
Let $A$ be the matrix of $p(t)$ whose its Jordan decomposition is given by
\[
    A \sim 
   \widetilde{A}= \begin{pmatrix}
        r_1 & 0 \\[4pt]
        0 & r_2
    \end{pmatrix},
\]
where $r_1$ and $r_2$ are in a quadratic extension $F$ of $K$. Define the spectral projections
\[
X = \frac{\widetilde{A} - r_2 I}{r_1 - r_2}, 
\qquad 
Y = \frac{\widetilde{A} - r_1 I}{r_2 - r_1}.
\]
It is easy to verify that
\[
X+Y=I, \quad X\widetilde{A} = r_1 X, \quad Y\widetilde{A} = r_2 Y, \quad YX=XY = 0, \quad X^2 = X,\; Y^2 = Y,
\]
and hence
\[
\widetilde{A}^n = r_1^n X + r_2^n Y, \qquad \text{for all } n \geq 0.
\]
Now, for any $r \in F$ we have the identity
\[
\prod_{\mu \in K}(t - \mu r) \;=\; t^q - r^{\,q-1} t,
\]
in the $K$-subspace $rK \simeq K$. Thus $t^q - r^{\,q-1} t$ vanishes on $\{\mu r \mid \mu \in K\}$. In particular, for $r = r_1$ or $r_2$, this polynomial annihilates all multiples $\lambda r_i$ with $\lambda \in K$. Since these polynomials do not necessarily have coefficients in $K$, we instead consider their product:
\[
z(t) = \bigl(t^q - r_1^{\,q-1} t\bigr)\bigl(t^q - r_2^{\,q-1} t\bigr)
= t^{2q} - (r_1^{\,q-1}+r_2^{\,q-1}) t^{q+1} + (r_1 r_2)^{q-1} t^2.
\]
As $r_1 r_2 = -\alpha$, and a computation shows
\[
r_1^{\,q-1} + r_2^{\,q-1} = -\frac{1+2\alpha}{\alpha} \;\in K,
\]
we have $z(t)$ has coefficients in $K$. We set
\[
\widetilde{E} := \widetilde{A}^{2q} - \bigl(r_1^{\,q-1}+r_2^{\,q-1}\bigr)\widetilde{A}^{q+1} + \widetilde{A}^2.
\]
Since $\widetilde{A}^q = I - \widetilde{A}$ (because 
\(\widetilde{A}^q = r_2 X + r_1 Y = (1-r_1)X + (1-r_2)Y = I - \widetilde{A}\)),
it follows that
\[
\widetilde{A}^{2q} = (I-\widetilde{A})^2 = I - 2\widetilde{A} + \widetilde{A}^2,
\qquad
\widetilde{A}^{q+1} = (I-\widetilde{A})\widetilde{A} = \widetilde{A} - \widetilde{A}^2.
\]
Substituting into $E$ yields
\[
\widetilde{E} = I + (-2-c)\widetilde{A} + (2+c)\widetilde{A}^2,
\qquad c := r_1^{\,q-1}+r_2^{\,q-1}.
\]
Since $\widetilde{A}^2 = \widetilde{A} + \alpha I$, this simplifies to
\[
\widetilde{E}= \bigl(1+\alpha(2+c)\bigr) I.
\]
But $2+c = -1/\alpha$, hence $\widetilde{E}=0$. The last equality yields that
\[
E := A^{2q} - \bigl(r_1^{\,q-1}+r_2^{\,q-1}\bigr)A^{q+1} + A^2=0.
\]
Since $r_1, r_2 \notin K$, we have that both are nonzero. Therefore $\alpha \neq 0$ and $A$ is invertible. Thus,
\[A^{2q-1} - \bigl(r_1^{\,q-1}+r_2^{\,q-1}\bigr)A^{q} + A=0\]
Then, for all $\lambda \in K$,
\begin{align*}
    e_1(\lambda e_3)^{(2q-1)} & - (r_1^{\,q-1}+r_2^{\,q-1}) e_1(\lambda e_3)^{(q)} + e_1(\lambda e_3)  \\
    & =\lambda \left(A^{2q-1} - \bigl(r_1^{\,q-1}+r_2^{\,q-1}\bigr)A^{q} + A\right) (e_1)= 0,
\end{align*}
and similarly with $e_2$ in place of $e_1$. This proves that the identity \eqref{RR7fin3} holds in $L_{10}$, since the polynomial is linear in $x_1$. 

This completes the proof.
\end{proof}

It is important to emphasize that the identity \eqref{RR7fin4} 
% and \eqref{RR7fin5}
, as well as \eqref{RR7fin3} 
% and \eqref{RR7fin6} 
(in their respective cases), play a fundamental role in bounding the exponents of the first and second variables in any monomial.
Indeed, modulo these identities, every monomial can be reduced so that the occurrences of the first two variables lie within a controlled range.

The following result formalizes this reduction property and constitutes a key step toward establishing the main theorem of this section.

\begin{lemma}\label{identi2}
Let $K$ be a finite field with cardinality $|K| = q$, and let $\alpha$ be the
parameter appearing in the multiplication table of the algebra $L_{10}$, so that
$r_1$ and $r_2$ are the roots of the polynomial $p(t)=t^{2}-t-\alpha$.

Assume either that:
\begin{itemize}
    \item both roots $r_1$ and $r_2$ do not lie in $K$, or
    \item the roots coincide in $K$ (in which case we additionally assume 
          $\operatorname{char}(K)\neq 2$).
\end{itemize}

Let $m_1,m_2,m_1',m_2'$ be positive integers bounded by $2q$ and $m_1\geq m_1^\prime$ such that
\[
    m_1 \equiv m_1' \pmod{q-1},
    \qquad
    m_2 \equiv m_2' \pmod{q-1}.
\]
Suppose there exist scalars $\lambda_1,\lambda_2\in K$ for which the polynomial
\[
    \lambda_1\, x_1^{(m_1)} x_2^{(m_2)} x_3^{(m_3)}\cdots x_n^{(m_n)}
    \;-\;
    \lambda_2\, x_1^{(m_1')} x_2^{(m_2')} x_3^{(m_3)}\cdots x_n^{(m_n)}
\]
is an identity of $L_{10}$. Then this polynomial belongs to the $T$-ideal generated by 
\eqref{RR7fin1}, \eqref{RR7fin2}, and \eqref{RR7fin4} (or by 
\eqref{RR7fin3}, depending on which of the above hypotheses on $r_1$ and $r_2$
holds).
\end{lemma}
\begin{proof}
Since we only need to analyze the first and second variables, it suffices to prove the statement for $n=2$.  

\medskip
\noindent
\textbf{Case 1:} $m_1' > 1$.  
In this case, the identity
\[
x_1^2 x_2^{(q)} - x_1^{(q+1)} x_2
\]
follows from \eqref{RR7fin2}, which implies $m_2 = m_2' < q$. Hence, we may assume $m_2 = m_2' = 1$ and consider the polynomial
\[
\lambda_1 x_1^{(m_1)} x_2 - \lambda_2 x_1^{(m_1')} x_2.
\]

Evaluating $x_1 \mapsto e_i + e_3$ for $i=1,2$ and $x_2 \mapsto e_3$, and using the notation of Lemma~\ref{startsame}, we obtain
\(
\lambda_1 x_1^{(m_1)} x_2 \mapsto \lambda_1 (e_i + e_3)^{(m_1)} e_3 = \lambda_1 A^{m_1}(e_i),
\)
and similarly
\[
\lambda_2 x_1^{(m_1')} x_2 \mapsto \lambda_2 A^{m_1'}(e_i).
\]
Thus, the identity assumption becomes
\[
A^{m_1'}\bigl(\lambda_1 A^{m_1-m_1'} - \lambda_2 I_2 \bigr)(e_i) = 0.
\]

Since $m_1' \ge 1$ and $A$ is invertible, we have
\(
(\lambda_1 A^{m_1-m_1'} - \lambda_2 I_2)(e_i) = 0.
\)
As $e_1$ and $e_2$ are linearly independent, it follows that
\[
\lambda_1 A^{m_1-m_1'} - \lambda_2 I_2 = 0.
\]
If $\lambda_1 \neq 0$, then
\(
A^{m_1-m_1'} = (\lambda_2 \lambda_1^{-1}) I_2.
\)
But since $m_1 \equiv m_1' \pmod{q-1}$ and $m_1 - m_1' < q$, we have $m_1 - m_1' = q-1$, so that
\[
(\lambda_2 \lambda_1^{-1}) A = A^{m_1-m_1'+1} = A^q = I - A.
\]
Applying this equality to $e_1$ gives 
\(
(\lambda_2 \lambda_1^{-1}) e_2 = e_1 - e_2,
\)
which is impossible because $e_1$ and $e_2$ are linearly independent. Therefore, we must have $\lambda_1 = \lambda_2 = 0$ in this case.

\medskip
\noindent
\textbf{Case 2:} $m_1' = 1$.  
Here, the evaluation $x_1 \mapsto e_1$ and $x_i \mapsto e_3$, for all $i>1$, forces $\lambda_1 = \lambda_2 = 0$, so the identity is trivial.

\medskip
Consequently, in all cases the polynomial lies in the $T$-ideal generated by the identities appearing in the hypothesis.
\end{proof}

We now have all the ingredients to complete the main objective of this subsection.
\begin{theorem}
Let $K$ be any field.
\begin{enumerate}
    \item If $K$ is infinite, then the $T$-ideal of all polynomial identities of $L_{10}$ is generated by the polynomial 
    \[
        x_1 (x_2 x_3).
    \]
    
    \item If $K$ is finite with cardinality $|K| = q$, and $\alpha \neq 0$ is chosen from the multiplication table of $L_{10}$ such that $r_1$ and $r_2$ are roots of $p(t) = t^2 - t - \alpha$, then the $T$-ideal $T(L_{10})$ is generated as follows:
    \begin{enumerate}
        \item[(i)] If $r_1$ and $r_2$ are distinct elements of $K$, by
        \[
            x_1 (x_2 x_3), \quad  \text{and} \quad x_1 x_2^{(q)} - x_1 x_2.
        \]
        
\item[(ii)] If $r_1 = r_2 = r$ in $K$ (i.e., $r$ is a repeated root), then 
$T(L_{10})$ is generated by the identities
\[
    x_1 (x_2 x_3), \qquad
    x_1 x_2 x_3^{(q)} - x_1 x_2^{(q)} x_3, \quad  \text{and} \quad 
    x_1 x_2^{(2q-1)} - 2\, x_1 x_2^{(q)} + x_1 x_2.
\]

        \item[(iii)] If both roots $r_1$ and $r_2$ do not lie in $K$, by
        \[
    x_1 (x_2 x_3), \qquad
    x_1 x_2 x_3^{(q)} - x_1 x_2^{(q)} x_3, \quad  \text{and} \quad 
     x_1 x_2^{(2q-1)} + \alpha^{-1} (1+2\alpha)\, x_1 x_2^{(q)} + x_1 x_2.
\]
    \end{enumerate}
\end{enumerate}
\end{theorem}
\begin{proof}
Firstly, Lemma \ref{startsame} implies that we may assume that any identity of $L_{10}$ begins with the variable $x_1$. 

If $K$ is an infinite field, then any multihomogeneous polynomial 
\(
f = f(x_1,x_2,\dots,x_n)
\) 
can be reduced, modulo $x_1(x_2 x_3)$, to a single monomial. Evaluating 
\[
x_1 \mapsto e_2 + e_3, \qquad x_k \mapsto e_3 \text{ for all } k \neq 1,
\] 
yields the desired result for the infinite field case.

Now assume that the field $K$ is finite, with cardinality $|K| = q$. Let $f$ be a regular identity of $L_{10}$.  
Without loss of generality, we can write
\[
f \equiv_{I_{L_{10}}} \sum_{\mathbf{m} = (\alpha_1, \alpha_2, \ldots, \alpha_n)} \lambda_{\mathbf{m}} \, x_1^{(\alpha_1)} x_2^{(\alpha_2)} \cdots x_n^{(\alpha_n)},
\]
where $\alpha_i < q$ for all $i > 2$ , by the identities \eqref{RR7cota} or \eqref{RR7fin2}. Here, $I_{L_{10}}$ denotes the $T$-ideal generated by the identities corresponding to each case.

If the polynomial \(p(t)=t^2-t-\alpha\) has distinct roots in \(K\), then \(\alpha_1\leq q\) and \(\alpha_i<q\) for all \(i>2\). The result therefore follows from the evaluation
\[
x_1 \mapsto e_2+\lambda_1e_3,
\qquad
x_i \mapsto \lambda_i e_3 \quad \text{for all } i>1.
\]

Now, assume that $r_1 = r_2 = r \in K$. We can see that the identity
$$x_1^{2q} - 2 x_1^{q+1} + x_1^{2}$$
is a consequence of \eqref{RR7fin4}. Together, these two identities imply that $\alpha_1 < 2q$ and $\alpha_2 < 2q-1$. Using the evaluations $x_1 \mapsto e_2$ \ (or $e_1$), $x_i \mapsto \lambda_i e_3$ for all $i > 1$, we see that any subpolynomial of the form
\[
\lambda_1 x_1 x_2^{(\alpha_2)} x_3^{(\alpha_3)} \cdots x_n^{(\alpha_n)}
-
\lambda_2 x_1 x_2^{(\alpha_2^{\prime})} x_3^{(\alpha_3)} \cdots x_n^{(\alpha_n)},
\]
for some $\lambda_1, \lambda_2 \in K$ and $\alpha_2 \equiv \alpha_2' \pmod{q-1}$, is an identity of $L_{10}$. Then, by Lemma \ref{identi2}, the desired result follows. It remains to consider the case where $x_1$ is not linear in $f$. From the identity
\[
x_1 x_2 x_3^{(q)} - x_1 x_2^{(q)} x_3,
\]
we see that 
\[
x_1^{(q+1)} x_2 - x_1^2 x_2^{(q)} \in T(L_{10}),
\] 
so we can assume $\alpha_2 = \alpha_2' < q$. Since $x_1$ is not linear in $f$, the conclusion follows from the evaluations
\[
x_1 \mapsto e_2 + \lambda_1 e_3, 
\qquad 
x_i \mapsto \lambda_i e_3 \quad \text{for all } i > 1,
\]
together with Lemma \ref{identi2}, and the fact that $1 < \alpha_1 < 2q$ for any $\alpha_1$ appearing in $\mathbf{m}$.

Finally, the case where $p(t)$ has no roots in $K$ is analogous to the preceding one, but relies on the identity \eqref{RR7fin3} instead of \eqref{RR7fin4}. Note that \[
r_1^{,q-1} + r_2^{,q-1} = -\frac{1+2\alpha}{\alpha} ;\in K.
\]

Thus, the desired result holds in all cases, which completes the proof.
\end{proof}

\begin{corollary}\label{BARL L_{10}}
Let $\mathcal{B}_{L_{10}}\subset \mathcal{L}\langle X\rangle$ be the set of all
monomials of the form
\[
    x_{j_1}^{(m_1)} x_{j_2}^{(m_2)} \cdots x_{j_n}^{(m_n)},
    \qquad
    (j_2 < j_3 < \dotsb < j_n 
    \ \text{and}\ 
    j_1 \neq j_k\ \text{for all}\ k\ge 2).
\]

If $K$ is an infinite field, then the images of $\mathcal{B}_{L_{10}}$ in the
relatively free algebra 
\(
\mathcal{L}\langle X\rangle / T(L_{10})
\)
form a basis of this algebra.

If $K$ is a finite field with $|K| = q$, we impose the constraints
$m_i < q$ for all $i = 3,\dots,n$, and additional restrictions depending on the
scalar $\alpha$:

\begin{enumerate}
    \item If the polynomial $p(t)=t^{2}-t-\alpha$ has
    distinct roots in $K$, then the images of $\mathcal{B}_{L_{10}}$ form a basis
    of the relatively free algebra provided that \(m_1 \le q\) and \(m_2 < q\).

    \item Otherwise (i.e., if $p(t)$ has no roots in $K$ or has a double root in $K$,
    with $\alpha\neq 0$), the following restrictions hold:
    \[
        m_2 < 2q - 1, \qquad \text{if } m_1 = 1,
    \]
    and, in all other cases, \(m_2 < q\) and \(1<m_1<2q\).
\end{enumerate}

Moreover, the codimension sequence of $L_{10}$ satisfies
\(
    c_n(L_{10})=n,
\)
and the exponent of the variety generated by $L_{10}$ is equal to~$1$.
\end{corollary}
\medskip

\subsection{\texorpdfstring{The algebra $L_{11}$}{The algebra L11}}\

Consider the multiplication in $L_{11}$ on the basis 
\(\{e_1,e_2,e_3\}\) given by  
\[
e_1 e_3 = e_1, 
\qquad  
e_2 e_3 = e_2.
\]

From a simple computation, we see that \(x_1 (x_2 x_3)\) is an identity for this algebra. Moreover, if $K$ is a finite field with $|K| = q$, \(x_1 x_2^{(q)} - x_1 x_2\) is also identity of $L_{11}$. Let $I_{L_{11}}$ the $T$-ideal generated by this identities. 

\begin{theorem}\label{theo:idRR11}
If $K$ is an infinite field, then the polynomial $x_1(x_2 x_3)$ generates the $T$-ideal of polynomial identities of $L_{11}$. If $|K| = q < \infty$, then one must add the polynomial $x_1 x_2^{(q)} - x_1 x_2$ to the generating set.
\end{theorem}
\begin{proof}
    Fixing an arbitrary $r \in \{1,\dots,n\}$ and setting $x_r \mapsto e_1 + e_3$ while keeping $x_i = e_3$ for $i \neq r$ gives $\lambda_r = 0$, and thus $f$ is a consequence of $x_1(x_2 x_3)$.

Now, suppose that $K$ is a finite field with $|K| = q$. For a regular identity $f(x_1,\dots,x_n)$ of $L_{11}$, we have
\[
f \equiv_{I_{L_{11}}} \sum_j \lambda_j x_{r_j}^{(m_{r_j j})} x_1^{(m_{1j})} \dotsb \widehat{x_{r_j}^{(m_{r_j j})}} \dotsb x_n^{(m_{n j})},
\]
where $m_{r_j j} \le q$ and $1 \le m_{i j} < q$ for all $i \ne r_j$.

For a fixed $k \in \{1,\dots,n\}$, let
\(
S_k = \{ j \mid r_j = k \}.
\)
We introduce new commutative variables $\beta_i$ and considering the evaluation
$x_k \mapsto e_1 + \beta_k e_3$ and $x_i \mapsto \beta_i e_3$ for  $i \ne k$. Under this evaluation, we obtain
\[
\sum_{j \in S_k} \lambda_j  \beta_k^{m_{k j} - 1} \beta_1^{m_{1 j}} \dotsb \widehat{\beta_k^{m_{k j}}} \dotsb \beta_n^{m_{n j}} = 0.
\]
From the known identities for $K$, it follows that $\lambda_j = 0$ for all $j \in S_k$. Since $k$ is arbitrary, we conclude that $f \in I_{L_{11}}$, which completes the proof.
\end{proof}

The algebra \(L_{11}\) contains the subalgebra $\langle e_1,e_1+e_3\rangle$, which is isomorphic to \(\mathcal{L}_4\); see \eqref{eq:L4}. Hence, $T(L_{11})\subseteq T(\mathcal{L}_4)$. On the other hand, the arguments used in the proof of Theorem~\ref{theo:idRR11} apply, with only minor modifications, to \(\mathcal{L}_4\) and show that no additional polynomial identities occur. Therefore,

$$
T(\mathcal{L}_4)=T(L_{11}).
$$

In particular, \(\mathcal{L}_4\) and \(L_{11}\) are PI-equivalent. Together with the previous results, this also yields the PI-equivalence of \(L_{11}\) and \(L_{10}(0)\).

\begin{corollary}
Over any field \(K\), the algebras \(L_{11}\), \(\mathcal{L}_4\), $N_4$ and \(L_{10}(0)\) are PI-equivalent. Consequently, their relatively free algebras coincide:

$$
\frac{\mathcal L\langle X\rangle}{T(L_{11})}
=
\frac{\mathcal L\langle X\rangle}{T(\mathcal L_4)}
=
\frac{\mathcal L\langle X\rangle}{T(N_4)}=
\frac{\mathcal L\langle X\rangle}{T(L_7(0))}.
$$
More explicitly, let \(\mathcal B_4\subseteq\mathcal L\langle X\rangle\) be the set of monomials of the form
$$
x_{j_1}^{(m_1)}x_{j_2}^{(m_2)}\cdots x_{j_n}^{(m_n)},
\qquad
j_2<j_3<\cdots<j_n,
\quad
j_1\neq j_k\ \text{for all }k\geq2.
$$
If \(K\) is infinite, then the images of the elements of \(\mathcal B_4\) form a basis of $\mathcal L\langle X\rangle/T(L_{11})$. If \(K\) is finite, with \(|K|=q\), the same conclusion holds upon imposing the additional restrictions
$$
1\leq m_1\leq q,
\qquad
1\leq m_i<q\quad (i\geq2).
$$
Moreover, $c_n(L_{11})=n$, and hence $\exp(L_{11})=1$. The same conclusions hold for $N_4$, \(\mathcal L_4\) and \(L_{10}(0)\).
\end{corollary}

\medskip

\subsection{\texorpdfstring{The algebra $L_{13}$}{The algebra L13}}\

In the algebra $L_{13}$, the multiplication on the basis $\{e_1,e_2,e_3\}$ is
\[
    e_3 e_3 = e_1,
    \qquad
    e_1 e_3 = e_1 + e_2 .
\]
A straightforward computation shows that, for arbitrary elements
\(
x_1=\alpha_1 e_1+\beta_1 e_2+\gamma_1 e_3
\)
and
\(
x_2=\alpha_2 e_1+\beta_2 e_2+\gamma_2 e_3,
\)
one has
\[
x_1 x_2
    = \gamma_2\bigl( (\alpha_1+\gamma_1)e_1 + \alpha_1 e_2 \bigr),
\]
and for any $s\ge 2$,
\[
x_1\, x_2^{(s)}
    = \gamma_2^{s} \bigl( (\alpha_1+\gamma_1)e_1 + (\alpha_1+\gamma_1)e_2 \bigr).
\]

Therefore, the polynomial
\(
    x_1 (x_2 x_3)
\)
is an identity of $L_{13}$.  
Moreover, when $|K| = q$, the polynomial
\(
    x_1 x_2^{(q)} - x_1 x_2 - x_2 x_1^{(q)} + x_2 x_1
\)
is also an identity of $L_{13}$.

Denote by $I_{L_{13}}$ the $T$-ideal generated by the above identities, according to the case under consideration.

Furthermore, the polynomials
\(
    x_1 x_2 x_3^{(q)} - x_1 x_2 x_3\) and \(
    x_1^{(q+1)} x_2 - x_1^{(2)} x_2\)
are in $I_{L_{13}}$.

\begin{lemma}\label{idg2RR9}
    In a finite field with cardinality $|K| = q$, every identity in two variables follows from the identities $x_1(x_2x_3)$ and $x_1 x_2^{(q)} - x_1 x_2 - x_2 x_1^{(q)} + x_2 x_1$.
\end{lemma}
\begin{proof}
For any regular identity \( f(x_1,\dots,x_n) \) of \(L_{13}\), we may write
\[
    f \equiv_{I_{L_{13}}} 
    \sum_{j \in J_1} \lambda_j^1\, x_{1}^{(m_{1j})} x_2^{(m_{2j})}
    \;+\;
    \sum_{k \in J_2} \lambda_k^2\, x_{2}^{(m_{1k})} x_1^{(m_{2k})},
\]
where, for every \(i \in J_1 \cup J_2\), the exponents satisfy  
\(1 \le m_{1i} \le q\) and \(1 \le m_{2i} \le q\).  
Moreover, if \(m_{1i} > 1\), then necessarily \(m_{2i} < q\); and in the special case  
\(m_{1i} = m_{2i} = 1\), we require that \(i \in J_1\).

Assume that there exists \(j_0 \in J_1\) such that \(m_{1j_0} = m_{2j_0} = 1\).
Under the evaluation
\(
x_1, x_2 \mapsto -e_1 + e_3,
\)
the corresponding monomial does not vanish, while all other monomials in the decomposition do vanish. Hence, \(\lambda^{1}_{j_0} = 0\).

Next, consider the evaluation
\(
x_1 \mapsto e_1, x_2 \mapsto \gamma_2 e_3, \gamma_2 \in K.
\)
This yields
\[
\sum_{j \in U_1} \lambda_j^1\, \gamma_2^{m_{2j}} = 0,
\]
where \(U_1 = \{ j \in J_1 \mid m_{1j} = 1 \}\).  
Since every \(j \in U_1\) satisfies \(2 \le m_{2j} \le q\), it follows that  
\(\lambda_j^1 = 0\) for all \(j \in U_1\).  
Hence, in the decomposition of \(f\), no monomial of the form  
\(x_1 x_2^{(m)}\) can appear. Now evaluate
\(
x_1 \mapsto (1 - \gamma_1)e_1 + \gamma_1 e_3, 
x_2 \mapsto -\gamma_2 e_1 + \gamma_2 e_3,
\gamma_1,\gamma_2 \in K.
\)
This produces the relation
\[
\sum_{j \in J_1} \lambda_j^1 \, 
    \gamma_1^{\,m_{1j}-1}\,
    \gamma_2^{\,m_{2j}}
    = 0.
\]
Since the exponents \(m_{1j}-1\) and \(m_{2j}\) range over a complete set allowed by the restrictions, we conclude that  
\(\lambda_j^1 = 0\) for all \(j \in J_1\).  
An identical argument shows that \(\lambda_k^2 = 0\) for all \(k \in J_2\). Therefore, all coefficients vanish and consequently \(f \in I_{L_{13}}\).
\end{proof}

\begin{theorem}
If $K$ is an infinite field, then the $T$-ideal of polynomial identities of the algebra $L_{13}$ is generated by the polynomial
\[
x_1(x_2 x_3).
\]

If $K$ is a finite field with $|K| = q$, then the $T$-ideal $T(L_{13})$ is generated by the polynomials
\[
x_1(x_2 x_3)
\qquad\text{and}\qquad
x_1 x_2^{(q)} - x_1 x_2 - x_2 x_1^{(q)} + x_2 x_1.
\]
\end{theorem}
\begin{proof}
Over an infinite field, we take a multihomogeneous identity 
\(f = f(x_1,x_2,\dots,x_n)\) of \(L_{13}\). 
As we have seen,
\[
    f \equiv_{I_{L_{13}}} 
    \sum_{j=1}^{n}\lambda_j\, x_j^{(m_j)} 
    x_1^{(m_1)} \cdots 
    \widehat{x_j^{(m_j)}} \cdots 
    x_n^{(m_n)}.
\]
If \(m_j = 1\) for some \(j \in \{1,\dotsc,n\}\), then evaluating 
\(x_j \mapsto e_1\) and \(x_k \mapsto e_3\) for all \(k \neq j\) gives 
\(\lambda_j = 0\). Thus, we may assume \(m_j > 1\) for all 
\(j = 1,\dotsc,n\).  

First, setting \(x_k \mapsto e_3\) for all \(k = 1,\dotsc,n\), we obtain  
\(\sum_{j=1}^{n}\lambda_j = 0\).  
On the other hand, evaluating 
\(x_{j_0} \mapsto e_1 + e_3\) and 
\(x_k \mapsto e_3\) for all \(k \neq j_0\), we obtain  
\(2\lambda_{j_0} + \sum_{j \neq j_0}\lambda_j = 0\).  
Hence \(\lambda_{j_0} = 0\), and since \(j_0\) is arbitrary, all 
coefficients vanish. This completes the argument over infinite fields.

It remains to consider the case where the field \(K\) is finite.
Assume that \(f(x_1,\dots,x_n)\) is a regular identity of \(L_{13}\).
By Lemma~\ref{idg2RR9}, we may suppose that \(n \ge 3\).
Then we can write
\[
f \equiv_{I_{L_{13}}} 
\sum_{j} \lambda_j\,
x_{r_j}^{(m_{r_j j})}
x_1^{(m_{1j})}
\cdots 
\widehat{x_{r_j}^{(m_{r_j j})}}
\cdots 
x_n^{(m_{nj})},
\]
where, for every \(j\), one has \(1 \le m_{r_j j} \le q\) and \(1 \le m_{ij} < q\) for all \(i \neq r_j\).

Fix \(k \in \{1,\dots,n\}\), and consider the set
$
J_k = \{\, j \mid r_j = k \,\}.
$
For arbitrary scalars \(\gamma_1,\dots,\gamma_n \in K\), evaluate the variables by
\(
x_k \mapsto (1 - \gamma_k)e_1 + \gamma_k e_3,
x_i \mapsto -\gamma_i e_1 + \gamma_i e_3
\ \ (i \ne k).
\)
Under this evaluation we obtain
\[
0 = 
\sum_{j \in J_k}
\lambda_j
\,
\gamma_k^{\,m_{kj} - 1}
\gamma_1^{m_{1j}}
\cdots
\widehat{\gamma_k^{m_{kj}}}
\cdots
\gamma_n^{m_{nj}}.
\]

Since the monomials in the scalars \(\gamma_1,\dots,\gamma_n\) appearing in the sum are distinct (by the bounds on the exponents), the above equality implies that
\[
\lambda_j = 0 \qquad\text{for all } j \in J_k.
\]

As this holds for every \(k \in \{1,\dots,n\}\), all coefficients \(\lambda_j\) vanish, and therefore \(f \in I_{L_{13}}\). 

This completes the proof.
\end{proof}

\begin{corollary}\label{BARL L_{13}}
Let $\mathcal{B}_{L_{13}} \subset \mathcal{L}\langle X\rangle$ be the set consisting of all monomials of the form
\[
x_{j_1}^{(m_1)} x_{j_2}^{(m_2)} \cdots x_{j_n}^{(m_n)}
\qquad\text{with}\qquad
j_2 < j_3 < \dotsb < j_n
\quad\text{and}\quad
j_1 \neq j_k \text{ for all } k \ge 2.
\]
If $K$ is an infinite field, then the images of $\mathcal{B}_{L_{13}}$ in the relatively free algebra 
$\mathcal{L}\langle X\rangle / T(L_{13})$ form a basis of this algebra.

If $K$ is a finite field with $|K| = q$, then the following additional restrictions on the exponents hold:
\begin{itemize}
    \item if $m_1 > 1$, then $m_2 < q$;
    \item if $n = 1$, then $1 \le m_1 \le q+1$;
    \item if $n = 2$, then $1 \le m_i \le q$ for $i = 1,2$, and in the case $m_1 = m_2 = 1$ we further require $j_1 < j_2$;
    \item if $n \ge 3$, then $1 \le m_1 \le q$ and $1 \le m_i < q$ for all $i \ge 2$.
\end{itemize}

Moreover, the codimension sequence of $L_{13}$ is given by
\(
c_n(L_{13}) = n,
\)
and the exponent of the variety $L_{13}$ is equal to \(1\).
\end{corollary}

	\bibliography{ref2}

  \bibliographystyle{amsplain.bst}
\end{document}